\documentclass[a4paper,12pt,twoside]{article}
\usepackage{pst-fill,pst-grad}
\usepackage{textcomp}
\usepackage[english]{babel}
\usepackage[utf8x]{inputenc}
\usepackage{graphicx}
\usepackage{amsmath}
\usepackage{multicol}
\usepackage{float}
\usepackage{fancyhdr}
\usepackage[matrix,arrow,curve]{xy}
\usepackage{pstricks} 
\usepackage{amsmath,amsfonts,verbatim,afterpage,theorem,euscript,mathrsfs,amssymb}
\usepackage{amsfonts}
\usepackage{amssymb}
\usepackage{array}
\usepackage{dsfont}
\usepackage{hyperref}
\newcommand{\mysection}{\setcounter{equation}{0} \section}

\newtheorem{Definition}{Definition}[section]
\newtheorem{Proposition}{Proposition}[section]
\newtheorem{Lemme}{Lemma}[section]
\newtheorem{Theoreme}{Theorem}

\newtheorem{Remarque}{Remark}[section]

\newcommand\R{\mathbb{R}}
\title{\bf Non Linear Singular Drifts and Fractional Operators: when Besov meets Morrey and Campanato}
\author{Diego Chamorro\footnote{Laboratoire de Math\'ematiques et Mod\'elisation d'Evry, Universit\'e d'Evry Val d'Essonne, France.}, St\'ephane Menozzi\footnote{Laboratoire de Math\'ematiques et Mod\'elisation d'Evry, Universit\'e d'Evry Val d'Essonne, France \& Laboratory of Stochastic Analysis, HSE, Moscow, Russia.}} 

\begin{document}
\maketitle
\begin{scriptsize}
\abstract{Within the global setting of singular drifts in Morrey-Campanato spaces presented in \cite{DCHSM}, we study now the H\"older regularity properties of the solutions of a transport-diffusion equation with nonlinear singular drifts that satisfy a Besov stability property. We will see how this Besov information is relevant and how it allows  to improve previous results. Moreover, in some particular cases we show that as the nonlinear drift becomes more regular, in the sense of Morrey-Campanato spaces, the additional Besov stability property will be less useful.}\\

\end{scriptsize}
\mysection{Introduction}
In a previous work \cite{DCHSM} we studied smoothness properties of the solutions $\theta:[0,+\infty[\times \mathbb{R}^{n}\longrightarrow \mathbb{R}$ of a \emph{linear} transport-diffusion equation of the following form
\begin{equation}\label{Equation00}
\partial_t \theta-\nabla\cdot(v\,\theta)+\mathcal{L}^{\alpha_{0}}\theta=0,
\end{equation}
with a divergence-free velocity field  $v:[0,+\infty[\times \mathbb{R}^{n}\longrightarrow \mathbb{R}^{n}$ for which we had a uniform control in the space $L^{\infty}\big([0,T],M^{q,a}(\mathbb{R}^{n})\big)$, where $M^{q,a}(\mathbb{R}^{n}) $ stands for the Morrey-Campanato space with parameters $q,a$ whose definition is recalled below, and with a \textit{diffusion} term that was given by a general non-degenerate \textit{stable-like} L\'evy operator $\mathcal{L}^{\alpha_{0}}$ of Blumenthal-Getoor index $0<\alpha_{0}<2$, \emph{i.e.} in a neighborhood of $0$ its symbol $\widehat {\mathcal L}^ {\alpha_0} $ is upper and lower bounded by $|\xi|^{\alpha_0} $. We will refer to $\alpha_0$ as the smoothness degree of the operator $\mathcal{L}^{\alpha_{0}} $.

We proved in this article that if we take as starting point an initial data $\theta_{0}\in L^{\infty}(\mathbb{R}^{n})$ and if the linear drift $v$ belongs to a Morrey-Campanato space $M^{q,a}({\mathbb R}^n)$ whose parameters $q,a$ and the smoothness degree $\alpha_{0}$ of the L\'evy operator are related with the dimension $n$ by the relationship
\begin{equation}\label{RelationMorrey}
\frac{a-n}{q}=1-\alpha_{0},
\end{equation}
then, weak solutions of equation (\ref{Equation00}) belong to a H\"older space $\mathcal{C}^{\gamma}(\mathbb{R}^{n})$ with $0<\gamma<1$. In the case of the fractional Laplacian, relation \eqref{RelationMorrey} readily follows from the homogeneity structure of equation \eqref{Equation00}. 

This actually yields a gain in regularity as $\theta_{0}\in L^{\infty}(\mathbb{R}^{n})$. Remark moreover that a small regularity index $\gamma$ is enough for our purposes since it is possible to apply bootstrap arguments to obtain higher regularity, see \cite{Caffarelli}, Section B, for the details. It is also worth noting here that the relationship (\ref{RelationMorrey}) is given by the equation's homogeneity and without any further information over the velocity field $v$ it is not possible to break it down as in some particular cases counterexamples can be exhibited, see \cite{Silv}.\\

We will adopt from now on the following notation: the index $\alpha_{0}$ will denote the smoothness degree of the operator $\mathcal{L}^{\alpha_{0}}$ in the \emph{linear} setting given by equation (\ref{Equation00}) and it will be always related to the dimension $n$ and the parameters of the Morrey-Campanato space $q,a$, by the equation (\ref{RelationMorrey}).\\

Recall now that L\'evy-type operators $\mathcal{L}^{\alpha_{0}}$ are, roughly speaking, generalizations of the fractional Laplacian $(-\Delta)^{\frac{\alpha_{0}}{2}}$ and the choice of Morrey-Campanato spaces $M^{q,a}$ (with $1<q<+\infty$ and $0\leq a<n+q$) was given by the fact that they can describe both singular and regular situations following the values of the parameters $q$ and $a$: indeed, if $0\leq a <n$ the corresponding spaces are singular (Morrey spaces) while if $n<a<n+q$ we obtain the following identification with H\"older spaces $M^{q,a}(\mathbb{R}^{n})\simeq \mathcal{C}^{\lambda}(\mathbb{R}^{n})$ where $\frac{a-n}{q}=\lambda$ and $\lambda\in ]0,1[$. The case $a=n$ is given by the space $bmo(\mathbb{R}^{n})$. See Section \ref{SeccDefinition_Resultats} below for a precise definition of these operators and spaces.\\

The aim of this article is then to consider singular and \emph{nonlinear} divergence-free velocity fields with some additional stability properties in order to break down the relationship (\ref{RelationMorrey}) in the following sense: we want to obtain a smaller regularity degree $0<\alpha<\alpha_{0}<2$ for the operator $\mathcal{L}^{\alpha}$ such that we still can obtain H\"older regularity for the associated solutions. The main ingredient to achieve this goal is to use the extra information given by the \textcolor{black}{boundedness} properties of the nonlinear drift (given in terms of Besov spaces): indeed, for a function $\theta:[0,+\infty[\times \mathbb{R}^{n}\longrightarrow \mathbb{R}$, with $n\geq 2$, we consider the problem
\begin{equation}\label{Equation0}
\partial_t \theta-\nabla\cdot(A_{[\theta]}\,\theta)+\mathcal {L}^{\alpha}\theta=0,
\end{equation}
where $\mathcal {L}^{\alpha}$ is a L\'evy-type operator of regularity degree $0<\alpha<2$ and $A_{[\theta]}$ is a nonlinear velocity field satisfying the following general hypotheses:
\begin{itemize}
\item $A_{[\theta]}$ is given as a vector of singular integrals of convolution type: 
\begin{equation}\label{DefinitionNoyau}
A_{[\theta]}(t,x)=[A_{1}(\theta)(t,x),\cdots, A_{n}(\theta)(t,x)],\qquad (n\geq 2),
\end{equation}
with $\displaystyle{A_{i}(\theta)(t,x)=v.p. \int_{\mathbb{R}^{n}}\kappa_{i}(x-y)\theta(t,y)dy}$  and where $\kappa_{i}=\kappa_{i}(x)$ is the associated kernel of the singular integral operator $A_{i}$ for $1\leq i\leq n$.
\item The nonlinear drift $A_{[\theta]}$ is divergence free.
\item The nonlinear drift $A_{[\theta]}$ satisfies some boundedness properties in Morrey-Campanato spaces $M^{q,a}({\mathbb R}^n)$ with $1\leq q<+\infty$ and $0\leq a<n$.
\item The operator $A_{[\theta]}$ is bounded in Besov spaces $\dot{B}^{s,p}_{q}({\mathbb R}^n)$ with $0\leq s<1$ and $1\leq p,q\leq +\infty$.
\end{itemize}
Let us make some remarks about the hypotheses over the nonlinear drift  $A_{[\theta]}$.  First, the singular integral setting given by (\ref{DefinitionNoyau}) allows us to consider a wide range of operators with remarkable properties studied in many books (see for example \cite{Stein2}, \cite{Torchinski}, \cite{Adams2}). Next, the divergence free condition is very natural in problems arising from fluid dynamics and this property is crucial in this work as it considerably simplifies the computations. \\

These types of nonlinear drifts appear in many applicative fields. We mention for instance the surface quasi-geostrophic equation (SQG), see e.g. \cite{Caffarelli}, for which $A_{[\theta]}(t,x)=[-R_2 \theta(t,x), R_1\theta(t,x)] $ where $(R_j)_{j\in \{1,2\}}$ denote the Riesz transforms \emph{i.e.} $\widehat {R_j \theta}(t,\xi)= -i \frac{\xi_j}{|\xi|}\widehat \theta(t,\xi)$. Other examples derived from magnetohydrodynamic equations (MHD) have been studied in the limit case $\alpha=2 $ by Friedlander and Vicol \cite{frie:vico:12}. They consider drifts of the form $(A_{[\theta]})_j(t,x)=\displaystyle{\sum_{i=1}^n}\partial_iT_{ij}\theta(t,x) $ where $ \{T_{ij}\}$ is an $n\times n $ matrix of Calder\'on-Zygmund singular integral operators such that $\displaystyle{\sum_{i,j=1}^n} \partial_{i}\partial_j T_{ij}=0 $, which guarantees the divergence free condition.\\

The boundedness of the drift in terms of Morrey-Campanato spaces was already used (in a linear setting) in \cite{DCHSM} but it will be used here in order to ensure existence issues (but not only!) and we will see how this boundedness hypothesis interacts with the solution's regularity.\\ 

The main novelty of this article relies in the handling of a Besov space boundedness condition for the nonlinear velocity field $A_{[\theta]}$ \textcolor{black}{combined with the maximum principle $\|\theta(\cdot, \cdot)\|_{L^{\infty}(L^{p})}\leq \|\theta_{0}\|_{L^{p}}$ which is satisfied by (weak) solutions of equation (\ref{Equation0})}.
Indeed, for a fixed final $T>0$ we can prove from the study of the maximum principle (see Section \ref{SecExiUnic} of the current work and Sections 2.2 and 2.3 in \cite{DCHSM} for details, \textcolor{black}{see also \cite{Cordoba}}) that
$$ \|\theta(T,\cdot)\|_{L^p}^p +p\int_0^T \int_{\mathbb{R}^n}|\theta(t,x)|^{p-2}\theta(t,x)\mathcal{L}^{\alpha} \theta(t,x)dxdt
\le \|\theta_0\|_{L^p}^p,$$
\textcolor{black}{from which we obtain the control
\begin{equation}\label{CTR_LP}
p\int_0^T \int_{\mathbb{R}^n}|\theta(t,x)|^{p-2}\theta(t,x)\mathcal{L}^{\alpha} \theta(t,x)dxdt
\le \|\theta_0\|_{L^p}^p.
\end{equation}
}
Recalling as well from \cite{DCHSM} \textcolor{black}{(see also \cite{PGDCH}) that for $0<t<T$ we have}
$$\|\theta(t,\cdot)\|_{\dot{B}^{\frac{\alpha}{p},p}_{p}}^p\le C\Big(\|\  |\theta(t,\cdot)|^{p/2}\ \|_{L^2}^2 +\int_{\mathbb{R}^n}|\theta(t,x)|^{p-2}\theta(t,x)\mathcal{L}^{\alpha} \theta(t,x)dx\Big),$$
\textcolor{black}{we thus derive from the maximum principle and from \eqref{CTR_LP} the following bound}
\begin{equation}
\label{CTR_BESOV_LP}
\begin{split}
\int_0^T \|\theta(t,\cdot)\|_{\dot{B}^{\frac{\alpha}{p},p}_{p}}^p dt 
&\le C(T)\|\theta_0\|_{L^p}^p,
\end{split}
\end{equation}
yielding $\theta \in L^p([0,T],\dot{B}^{\frac{\alpha}{p},p}_{p}) $. 
The fact that the quantity in the l.h.s. of \eqref{CTR_LP} 
 expresses regularity in terms of Besov spaces was already observed in \cite{CW} (see as well \cite{PGDCH}, \cite{DCHSM} and Theorem \ref{RegulariteBesov} below) but it seems that it was not fully used to study the regularity of the solution of this type of equations. We mention that the regularity conveyed by the Besov space $\dot{B}^{\frac{\alpha}{p},p}_{p}$ is only related to the initial data $\theta_{0}\in L^{p}$ and to the smoothness degree $\alpha$ of the L\'evy-type operator. \textcolor{black}{Now, in order to link this Besov information with the behavior of the drift, we} will assume the following boundedness condition:
 \begin{equation}\label{Control_Besov1}
\|A_{[\theta]}\|_{L^p([0,T], \dot{B}^{s,p}_{q})}\leq C\|\theta\|_{L^p([0,T],\dot{B}^{\frac \alpha p,p}_{p})}.
\end{equation}
In this article we will show how to use this additional Besov regularity given by (\ref{Control_Besov1}) and 
\eqref{CTR_BESOV_LP}
in order to by-pass the Morrey-Campanato boundedness hypothesis and we will study how to break down the relationship (\ref{RelationMorrey}): within this general framework we will see that if the initial data $\theta_{0}$ belongs to a Lebesgue space $L^{p}$ and if the parameter $p$ that characterizes this initial condition, the Morrey-Campanato indexes $q,a$ and the regularity degree $\alpha$ of the L\'evy-type operator are \emph{suitably related}, then it is possible to obtain H\"older regularity for the solutions of equation (\ref{Equation0}) and moreover we obtain that the corresponding value of the regularity degree $\alpha$ satisfies $\alpha<\alpha_{0}$ where $\alpha_{0}$ was given in relationship (\ref{RelationMorrey}): we will thus obtain results that are out of reach in our previous work \cite{DCHSM}.\\

The main ideas for the proofs are again based on the Hardy-H\"older duality which was already the cornerstone in \cite{DCHSM} and we now specifically exploit the additional Besov information to improve the results (see Remark \ref{dicho_MC_B} for details).\\

It is worth noting here that there are two situations where this extra regularity information is not as relevant as expected. The first one is related to inequalities (\ref{Control_Besov1}) and \textcolor{black}{(\ref{CTR_BESOV_LP})} which express the fact that the Besov regularity $\dot{B}^{\frac{\alpha}{p},p}_{p}$ is linked to the initial data $\theta_{0}\in L^{p}$ and if we made $p\longrightarrow +\infty$ we would completely loose this regularity information. 
The second one is related to the nonlinear drift: we are able to use this extra Besov regularity information as long as the drift belongs to a singular Morrey-Campanato space (\emph{i.e.} if $0\leq a<n$) but we do not know if it is possible to exploit the Besov information in a more regular Morrey-Campanato framework (\emph{i.e.} if $n< a<n+q$). Recall that when $n< a<n+q$, Morrey-Campanato spaces are equivalent to classical H\"older spaces since $M^{q,a}(\mathbb{R}^{n})\simeq \mathcal{C}^{\lambda}(\mathbb{R}^{n})$ where $\frac{a-n}{q}=\lambda$ with $\lambda\in ]0,1[$, and the Besov regularity is probably too weak to improve this already regular setting.\\

Our last remark concerns the fact that Besov and Morrey-Campanato boundedness properties are not very restrictive conditions as all standard (and reasonable) Calder\'on-Zygmund operators are bounded in these functional spaces. See \cite{Stein2} and \cite{Adams2} for more details about this fact. \\

The plan of the article is the following. In Section \ref{SeccDefinition_Resultats} we give the precise definition of L\'evy-type operators used in equation (\ref{Equation0}) and we recall the definition and some useful properties of Morrey-Campanato and Besov spaces. We also state in this section the main theorems and results of the article 
that emphasize that,
according to the assumptions on the nonlinear velocity field, different results are obtained. In particular, we will show then an interesting competition between the Morrey-Campanato information and the Besov regularity. In Section \ref{SecExiUnic} we discuss some existence issues for this nonlinear PDE. Finally, Section \ref{SeccHolderRegularity} is devoted to the proofs. 
\mysection{Notation and Presentation of the results}\label{SeccDefinition_Resultats}

We start with the definition of the L\'evy-type operator $\mathcal{L}^{\alpha}$ used in equation (\ref{Equation0}). These operators are natural generalizations of the fractional Laplacian $(-\Delta)^{\frac{\alpha}{2}}$, with $0<\alpha<2$, which is given in the Fourier variable by the expression $\widehat{(-\Delta)^{\frac{\alpha}{2}} \varphi}(\xi)=c|\xi|^{\alpha}\widehat{\varphi}(\xi)$, for all suitable regular enough functions $\varphi$. Thus, following this approach, we define the operator $\mathcal{L}^{\alpha}$ by the expression
\begin{equation}\label{DefinitionOperateur}
\widehat{\mathcal{L}^{\alpha}\varphi}(\xi)=a(\xi)\widehat{\varphi}(\xi), 
\end{equation}
here the symbol $a$ is given by the so-called Lévy-Khinchin formula 
$$\displaystyle{a(\xi)=\int_{\mathbb{R}^n\setminus \{0\}}\big(1-\cos(\xi\cdot y)\big)\pi(y)dy},$$
where the function $\pi$ is symmetric, i.e. $\pi(y)=\pi(-y) $ for all $y\in \R^n$ and satisfies the bounds
\begin{equation}\label{DefKernelLevy}
\begin{split}
\overline{c}_1|y|^{-n-\alpha}\leq &\pi(y)\leq \overline{c}_2|y|^{-n-\alpha} \qquad \mbox{over } |y|\leq 1,\\
0\leq &\pi(y)\leq \overline{c}_2|y|^{-n-\delta} \qquad \mbox{over } |y|> 1,
\end{split}
\end{equation}
here $0<\overline{c}_1\leq\overline{c}_2$ are positive constants and where $0<\delta <\alpha<2$. The important point here is that the parameter $\alpha$ (called the \emph{smoothness degree}) will rule the smoothness properties of the operator $\mathcal{L}^{\alpha}$ in the sense that for such $\alpha$ the regularizing effect of $\mathcal{L}^{\alpha}$ is similar to the fractional Laplacian $(-\Delta)^{\frac{\alpha}{2}}$. See \cite{Jacob} for additional properties concerning Lévy operators and the Lévy-Khinchin formula. \\

We define now local Morrey-Campanato spaces $M^{q,a}(\mathbb{R}^{n})$  for $1\leq q<+\infty$,  $0\leq a<n+q$  by the condition:
$$\|f\|_{M^{q,a}}=\underset{\underset{0<r<1}{x_0\in \mathbb{R}^n}}{\sup}\; \left(\frac{1}{r^{a}}\int_{B(x_0,r )}|f(x)-\overline{f}_{B(x_0,r)}|^{q} dx\right)^{1/q}+ \underset{\underset{r\geq 1}{x_0\in \mathbb{R}^n}}{\sup}\; \left(\frac{1}{r^{a}}\int_{B(x_0,r )}|f(x)|^{q}dx\right)^{1/q}<+\infty,$$
where $\overline{f}_{B(x_0,r)}=\displaystyle{\frac{1}{|B(x_{0},r)|}\int_{B(x_{0},r)}f(y)dy}$ is the average of the function $f$ over the ball $B(x_{0},r)$.\\

Remark that if $a=0$ and $1\leq q<+\infty$ we have $L^{q}\subset M^{q,0}$. When $a=n$ and $q=1$ the space $M^{1,n}(\mathbb{R}^{n})$ corresponds to the space $bmo$ (the local version of $BMO$) and from this fact we derive the identification $M^{1,n}\simeq M^{q,n}$ for $1< q<+\infty$.  Finally, if $n<a<n+q$ and $1\leq q<+\infty$, we obtain $M^{q,a}(\mathbb{R}^n)= \mathcal{C}^{\lambda}(\mathbb{R}^n)$,  classical H\"older space with $0<\lambda=\frac{a-n}{q}<1$. The case $n+q\leq a$ will not be considered as the corresponding spaces are reduced to constants.
\begin{Remarque}\label{Rem_MorreyRegularity}
 As we can see, if we consider the range of values $0\leq a<n$, the Morrey-Campanato space $M^{q,a}$ becomes \emph{less singular} as the value of the parameter $a$ increases.
\end{Remarque} 
Finally, when $0\leq a <n$, the H\"older inequality easily yields the following embedding :
\begin{equation}\label{Injection_MorreyLebesgue}
L^{p}(\mathbb{R}^{n})\subset M^{q,a}(\mathbb{R}^{n})\qquad \mbox{with } p=\frac{qn}{n-a}.
\end{equation}
This particular fact suitably combined with a maximum principle will be quite useful in the following lines. 
See \cite{Adams}, \cite{Peetre} and  \cite{Zorko} for more details concerning Morrey-Campanato spaces.\\

For $0<s< 2$ and $1\leq q \leq +\infty$, homogeneous Besov spaces $\dot{B}^{s,p}_{q}(\R^n)$ may be defined by the condition
$$\|f\|_{\dot{B}^{s,p}_{q}}=\left(\int_{\mathbb{R}^{n}}\frac{\|f(\cdot+y)+f(\cdot-y)-2f(\cdot)\|^{q}_{L^{p}}}{|y|^{n+sq}}dy\right)^{1/q}<+\infty.$$
However, in the particular case when $0<s<1$ and $p=q$, the previous condition can be replaced by the following one which will be used in the sequel:
\begin{equation}\label{DefinitionBesov}
\|f\|_{\dot{B}^{s,p}_{p}}=\left(\int_{\mathbb{R}^{n}}\int_{\mathbb{R}^{n}}\frac{|f(x)-f(y)|^{p}}{|x-y|^{n+sp}}dxdy\right)^{1/p}<+\infty.
\end{equation}
For more details about Besov spaces see the monograph \cite{RS}. \\

From now on we will always assume that $\mathcal{L}^{\alpha}$ is a L\'evy-type operator of the form (\ref{DefinitionOperateur}) satisfying hypothesis (\ref{DefKernelLevy}) with a smoothness degree $0<\alpha <2$ and for the nonlinear drift $A_{[\theta]}$ we assume that it is given as a divergence free vector of singular integrals of convolution type as in (\ref{DefinitionNoyau}). \\

With these ingredients, and for a initial data $\theta_{0}$, we consider now the following nonlinear equation for a function  $\theta:[0,+\infty[\times \mathbb{R}^{n}\longrightarrow \mathbb{R}$ with $n\geq 2$:
\begin{equation}\label{EquationUNA}
\begin{cases}
\partial_t \theta(t,x)-\nabla\cdot(A_{[\theta]}\,\theta)(t,x)+\mathcal{L}^{\alpha}\theta(t,x)=0,\\[4mm]
\theta(0,x)=\theta_{0}(x), \quad \mbox{for all } x\in \mathbb{R}^{n}.
\end{cases}
\end{equation}
\textcolor{black}{In the following, for a time-space function $v:[0,T]\times \R^n\rightarrow \R$, such that $v\in L^p([0,T],X) $, where $X$ is a function space and $p\in [1,+\infty] $, we write $\|v\|_{L^p(X)}:=\|v\|_{L^p([0,T],X)} $. We will take as possible spaces for $X$, $L^p=L^p(\R^n),\ M^{q,a}=M^{q,a}(\R^n),\ \dot{B}^{s,p}_{q}=\dot{B}^{s,p}_{q}(\R^n)$.}\\

We have now all the notation and definition that will allow us to state our first theorem.\\

\begin{Theoreme}[General Framework
]\label{Theoreme1} Consider equation (\ref{EquationUNA}) above and assume the following points:
\begin{itemize}
\item[1)] The initial data $\theta_{0}$ belongs to a Lebesgue space $L^{p}\cap L^{\infty}(\mathbb{R}^{n})$ with $n(n-1)\vee 2n<p<+\infty$,
\item[2)] The operator $A_{[\theta]}$ that defines the nonlinear drift in equation (\ref{EquationUNA}) above is
\begin{itemize}
\item[$(a_{1})$] bounded in Morrey-Campanato spaces $M^{q,a}$ 
with $(n-a)(n-1)\vee 2n \leq q<+\infty$ and $0\leq a <n$,
\begin{equation}\label{HypoMorreyBorne1}
\|A_{[\theta]}\|_{L^{\infty}(M^{q,a})}\leq C_{A}\|\theta\|_{L^{\infty}(M^{q,a})}\qquad (0\leq a <n),
\end{equation}
\item[$(b_{1})$] bounded in Besov spaces $\dot{B}^{\frac{\alpha}{p},p}_{p}
$,
\begin{equation}\label{HypoBesovBorne1}
\|A_{[\theta]}\|_{L^{p}(\dot{B}^{\frac{\alpha}{p},p}_{p})}\leq D_{A}\|\theta\|_{L^{p}(\dot{B}^{\frac{\alpha}{p},p}_{p})},
\end{equation}
where $1<\alpha<2$ is the regularity degree of the operator $\mathcal{L}^{\alpha}$ and $1<p<+\infty$ is given in the previous point 1) above.\\
\end{itemize}
\end{itemize}
If the index $p$ defining the initial data Lebesgue space $L^{p}$, the regularity degree $\alpha$ of the L\'evy-type operator  $\mathcal{L}^{\alpha}$ and the parameters $q,a$ of the Morrey-Campanato space $M^{q,a}$ are related by the conditions 
\begin{equation}\label{Equation_Existence}
p=\frac{qn}{n-a},\qquad \mbox{and}\qquad 1<\alpha=\frac{p+n}{p+1}<2,
\end{equation}
then the solution $\theta(t,\cdot)$ of the equation (\ref{EquationUNA}) belongs to a H\"older space $\mathcal{C}^{\gamma}(\mathbb{R}^{n})$ for $0<\gamma<\frac{1}{4}$ small. \\

Moreover, the regularity degree $\alpha$ given in (\ref{Equation_Existence}) using the Besov stability property of the drift is always less than the degree $\alpha_{0}$ given in (\ref{RelationMorrey}) which can be deduced from the Morrey-Campanato boundedness property.
\end{Theoreme}
Remark now that we take as initial data a function that belongs to the space $L^{p}\cap L^{\infty}$: the $L^{\infty}$ information is needed to apply the regularity results of \cite{DCHSM} and to obtain a general existence framework while the $L^{p}$ control is needed in order to use the Morrey-Campanato and Besov regularity information. \\

The lower bound $n(n-1)\vee 2n<p$ stated for the initial data is given in this form for the sake of simplicity and it is mostly related to computational issues. \textcolor{black}{It guarantees that the results from \cite{DCHSM} apply and that the Morrey-Campanato information (\ref{HypoMorreyBorne1}) yields to the relation \eqref{RelationMorrey} to derive the expected regularization property. A more accurate bound on $p$ would} depend on the smoothness degree $\alpha$ and some other technical parameters (see Section \ref{SeccHolderRegularity} and Theorem 2 of \cite{DCHSM})\footnote{\textcolor{black}{For Theorem 2 in \cite{DCHSM} to apply, the constraint $ \frac{\alpha}{\alpha-1}>\textcolor{black}{n}$ is needed. Since $\alpha=\frac{p+n}{p+1} $, this first gives $p>n(n-2) $. On the other hand, to compare with the usual homogeneity relation \eqref{RelationMorrey}, we need this constraint to be fulfilled as well by $ \alpha_0=1-\frac{a-n}{q}$. Since $p$ and $q$ are related by \eqref{Equation_Existence}, we then get $p>n(n-1)$. Eventually, the condition $q>\frac{n}{1-\gamma}$ also appears in Theorem 2 of \cite{DCHSM}. Choosing $\gamma<\frac 14 $, \eqref{Equation_Existence} yields that this last condition is satisfied for $p>2n $.}}.
Remark also that, since $p$ and $q$ are related by the first equality of (\ref{Equation_Existence}), the condition $(n-a)(n-1)\vee 2n \leq q<+\infty$ is a consequence of the previous bound on $p$.\\

Observe that we have imposed two conditions over the nonlinear drift: a Morrey-Campanato control (\ref{HypoMorreyBorne1}) and a Besov control (\ref{HypoBesovBorne1}) and it is interesting to study which one of theses conditions gives a better result in the sense that the corresponding smoothness degree $\alpha$ needed to deduce H\"older continuity is as low as possible.\\ 

From the Morrey-Campanato control (\ref{HypoMorreyBorne1}) and recalling that we have the first relationship of (\ref{Equation_Existence}) between the parameters $p,n,q$ and $a$, we can use the Morrey-Campanato-Lebesgue embedding (\ref{Injection_MorreyLebesgue}). Combining this fact with the maximum principle (see point $(iii)$, page \pageref{MaxPrinciple} and the computations performed in (\ref{PointFixe2})-(\ref{Formula_FormaBilineal})) we have $\|A_{[\theta]}\|_{L^{\infty}(M^{q,a})}\leq C_{A}\|\theta\|_{L^{\infty}(M^{q,a})}\leq C_{A}\|\theta\|_{L^{\infty}(L^{p})}\leq C_{A}\|\theta_{0}\|_{L^{p}}<+\infty$. Thus, we can apply the results given in \cite{DCHSM} to obtain a smoothness degree $\alpha_{0}$ given by the initial relationship (\ref{RelationMorrey}).\\

But using the Besov stability hypothesis (\ref{HypoBesovBorne1}) we have a stronger result since in this case we always obtain a corresponding smoothness degree $\alpha$ such that $1<\alpha<\alpha_{0}<2$. Indeed, following (\ref{RelationMorrey}) we have $\alpha_{0}=\frac{q+n-a}{q}$ but since $p$ is related to $q$ by (\ref{Equation_Existence}), we obtain that $\alpha_{0}=\frac{p+n}{p}>\frac{p+n}{p+1}=\alpha$: the Besov regularity information allows us to improve the general result.\\ 

Remark now that, by definition we have $\alpha=\frac{p+n}{p+1}>1$ and, although we can break down the relationship (\ref{RelationMorrey}), the effect of the Besov regularity information is not strong enough to obtain a corresponding smoothness degree $\alpha$  such that $0<\alpha<1$.\\

Note also that, if we let $a\longrightarrow n$, then the Morrey-Campanato space $M^{q,a}$ becomes \emph{less singular} and gets closer to the bounded mean oscillation space $bmo\simeq M^{q,n}$ (recall Remark \ref{Rem_MorreyRegularity}), we thus have from the relationship (\ref{Equation_Existence}) that $p\longrightarrow +\infty$ so that the additional regularity information given by (\ref{HypoBesovBorne1}) vanishes -we loose all the extra regularity information- and we obtain in this limit $\alpha \longrightarrow 1$. This framework corresponds to an $L^{\infty}$ initial data, a drift term in $bmo$ and a diffusion term of regularity degree $\alpha=1$ and this particular situation was first studied in \cite{Caffarelli} and then generalized in \cite{Chamorro}.\\

Theorem \ref{Theoreme1} gives a natural and general framework since conditions (\ref{HypoMorreyBorne1}) and (\ref{HypoBesovBorne1}) are not very restrictive ones. As said in the introduction, a wide family of Calder\'on-Zygmund operators satisfy these hypotheses (for more details see \cite{Adams2} for Morrey-Campanato spaces and \cite{PGLM1}, \cite{Meyer} for Besov spaces).\\

In the next theorem we will show an interesting competition between the Morrey-Campanato and the Besov regularity. Indeed, we will assume here a different Morrey-Campanato boundedness behavior for the nonlinear drift $A_{[\theta]}$ and following Remark \ref{Rem_MorreyRegularity} we will see that, under special hypotheses, as long as the parameter $a$ of the Morrey-Campanato space $M^{q,a}$ remains small enough (the Morrey-Campanato space is then more singular) the Besov regularity information is relevant, but as the parameter $a$ increases (the Morrey-Campanato space becomes less singular) the Besov regularity information will be useless.

\begin{Theoreme}[Morrey-Campanato and Besov competition]\label{Theoreme2} Let $\theta_{0}\in L^{p}\cap L^{\infty}(\mathbb{R}^{n})$, with $n\geq 2$, be an initial data with $n(n-1)\vee 2n<p<+\infty$ and consider the nonlinear equation (\ref{EquationUNA}). Assume the following points
\begin{itemize}
\item[1)] $\mathcal{L}^{\alpha}$ is a L\'evy-type operator of the form (\ref{DefinitionOperateur}) satisfying hypothesis (\ref{DefKernelLevy}) with a smoothness degree $\alpha$ such that
\begin{equation*}
1<\alpha=\frac{p+n}{p+1}<2.
\end{equation*} 
\item[2)] The operator $A_{[\theta]}$ that defines the nonlinear drift in equation (\ref{Equation01}) above is
\begin{itemize}
\item[$(a_{2})$] bounded in $L^{p}-M^{p,a}$ spaces in the following sense
\begin{equation}\label{HypoMorreyBorne2}
\|A_{[\theta]}\|_{L^{\infty}(M^{p,a})}\leq C_{A}\|\theta\|_{L^{\infty}(L^{p})}\qquad (
0\leq a <n), 
\end{equation} 
\item[$(b_{2})$] bounded in Besov spaces 
$$\|A_{[\theta]}\|_{L^{p}(\dot{B}^{\frac{\alpha}{p},p}_{p})}\leq D_{A}\|\theta\|_{L^{p}(\dot{B}^{\frac{\alpha}{p},p}_{p})}.
$$
\end{itemize}
\end{itemize}
Then, as long as we have the following condition between the Morrey-Campanato parameter $a$, the regularity degree $\alpha$ and the Lebesgue index $p$ of the initial data:
\begin{equation}\label{MorreyvsBesov1}
0\leq a< \alpha=\frac{p+n}{p+1},
\end{equation}
the solution $\theta(t,\cdot)$  of the equation (\ref{EquationUNA})  belongs to a H\"older space $\mathcal{C}^{\gamma}(\mathbb{R}^{n})$ for $0<\gamma<\frac{1}{4}$ small and moreover this value of $\alpha$ satisfies $1<\alpha<\alpha_{0}<2$ and we break down the relationship (\ref{RelationMorrey}).\\

However, if we have 
$$ a\geq \alpha=\frac{p+n}{p+1},$$ 
then the solution $\theta(t,\cdot)$  of the equation (\ref{EquationUNA}) still belongs to a H\"older space $\mathcal{C}^{\gamma}(\mathbb{R}^{n})$ for $0<\gamma<\frac{1}{4}$ small, but we have now $1<\alpha_{0}\leq \alpha<2$: the Besov regularity information is useless and we do not break down the relationship (\ref{RelationMorrey}).
\end{Theoreme}
We make here some remarks and comments. 
\begin{itemize}
\item[$\bullet$] As for Theorem \ref{Theoreme1},  the lower bound $n(n-1)\vee 2n<p$ is technical.
\item[$\bullet$] If $0\leq a < \alpha=\frac{p+n}{p+1}$, then it is easy to see that the parameters $a,q,n$ and $\alpha$ do not fulfill equation (\ref{RelationMorrey}) and that we have $1<\alpha<\alpha_{0}<2$. 
\item[$\bullet$] Morrey-Campanato spaces $M^{p,a}$ become more and more regular as the value of the parameter $a$ grows (recall Remark \ref{Rem_MorreyRegularity}). The smallness condition $0\leq a<\alpha$ stated in Theorem \ref{Theoreme2} reflects the fact that the Besov stability will help us only in very singular situations.
\item[$\bullet$] In particular, due to homogeneity arguments, if $\alpha< a<n+p$, the Besov stability is useless if we compare it to the Morrey-Campanato norm control given by the hypothesis (\ref{HypoMorreyBorne2}). Indeed, in this case the smaller smoothness degree for L\'evy operators will be given by the relationship (\ref{RelationMorrey}) instead of (\ref{MorreyvsBesov1}): this fact shows that condition (\ref{HypoMorreyBorne2}) implies regularization properties for the nonlinear drift $A_{[\theta]}$. See \cite{Ros} and the reference there in for a characterization of operators fulfilling this type of boundedness conditions. 
\item[$\bullet$] In the limit case when $a=\alpha=\frac{p+n}{p+1}$, then a Morrey-Campanato control and Besov stability provides the same regularity information: Besov finally meets Morrey-Campanato. \\[2mm]
\end{itemize}
As we can observe from Theorems \ref{Theoreme1} and \ref{Theoreme2}, the additional Besov regularity information allows us to break down to some extent the relationship (\ref{RelationMorrey}) and we can consider lower regularity degrees for the L\'evy-type operators, but only in the range $1<\alpha<2$, and this extra regularity is not strong enough to obtain results in the range $0<\alpha\leq 1$.\\

This last case is particularly of interest when investigating the supercritical SQG equation, i.e. for $0<\alpha<1 $. For the critical value $\alpha=1 $, from the seminal work of Caffarelli and Vasseur \cite{Caffarelli} we recall that if the initial data belongs to $L^{2}$, an $ L^\infty$ control holds for the solution. It is then known that the drift based on the Riesz transform belongs to $BMO$. This allows to derive, either by De Giorgi like techniques in \cite{Caffarelli} or through the molecular like decomposition, see e.g. \cite{KN}, \cite{Chamorro}, the H\"older continuity of the solution. However, when $0<\alpha<1$ by the homogeneity condition \eqref{RelationMorrey}, to derive the H\"older continuity of the solution would require to establish some ($1-\alpha$)-H\"older regularity properties on the drift, see \cite{CW} and \cite{CW2}. This is a very challenging open problem as we do not know if the drift of the SQG equation given in terms of Riesz transforms can be related (due to some cancellation property or to an extra hidden regularity property) to any kind of regularity. Therefore, a natural approach consists in considering a regularized drift.  This has been done for instance in \cite{Delga} for the linear case. We consider  below a regularized drift and discuss how to use the Besov information to break the expected homogeneity. As explained in the introduction, our approach has some interest if  $\theta_0\in L^p, p<+\infty $.

In order to investigate this situation and to obtain a better understanding of the Morrey-Campanato/Besov competition we state now our two last results which are variations of Theorem \ref{Theoreme1} and Theorem \ref{Theoreme2}. We will first introduce in Theorem \ref{Theoreme3} below a regularization in the nonlinear drift by considering $A_{[(-\Delta)^{\frac{-\eta}{2}}\theta]}$ where $0<\eta<1$ is a small fixed parameter.
\begin{Theoreme}[Regular Drift]\label{Theoreme3} Let $\theta_{0}\in L^{p}\cap L^{\infty}(\mathbb{R}^{n})$ be an initial data with $n(n-1)\vee 2n<p<+\infty$, let $0<\eta<1$ be a fixed (small) parameter and consider the nonlinear equation
\begin{equation}\label{Equation01}
\begin{cases}
\partial_t \theta(t,x)-\nabla\cdot(A_{[(-\Delta)^{-\frac{\eta}{2}}\theta]}\,\theta)(t,x)+\mathcal{L}^{\alpha}\theta(t,x)=0,\\[4mm]
\theta(0,x)=\theta_{0}(x). 
\end{cases}
\end{equation}
Assume the following points
\begin{itemize}
\item[1)] $\mathcal{L}^{\alpha}$ is a L\'evy-type operator of the form (\ref{DefinitionOperateur}) satisfying hypothesis (\ref{DefKernelLevy}) with a smoothness degree $\alpha$ such that
\begin{equation}\label{DegreLevy11}
1/2<\alpha=\frac{n+p(1-\eta)}{p+1}<2.
\end{equation} 
\item[2)] The operator $A_{[(-\Delta)^{-\frac{\eta}{2}}\theta]}$ that defines the nonlinear drift in equation (\ref{Equation01}) above 
\begin{itemize}
\item[$(a_{3})$] is bounded in Morrey-Campanato spaces 
\begin{equation}\label{HypoMorreyBorne3}
\|A_{[(-\Delta)^{-\frac{\eta}{2}}\theta]}\|_{L^{\infty}(M^{q,a})}\leq C_{A}\|(-\Delta)^{-\frac{\eta}{2}}\theta\|_{L^{\infty}(M^{q,a})}\qquad (0\leq a <n),
\end{equation} 
\item[$(b_{3})$] is bounded in Besov spaces $\|A_{[(-\Delta)^{-\frac{\eta}{2}}\theta]}\|_{L^{p}(\dot{B}^{\alpha/p\textcolor{black}{+\eta},p}_{p})}\leq D_{A}\|(-\Delta)^{-\frac{\eta}{2}}\theta\|_{L^{p}(\dot{B}^{\alpha/p\textcolor{black}{+\eta},p}_{p})}$\\ with $0< \alpha <2$.
\end{itemize}
\end{itemize}
Then if $0<\frac{\alpha}{p}+\eta<1$ and if the parameters $p$, $q\ge 2n,a$ and $\eta$ satisfy the identity
\begin{equation}\label{RelationMorreySobolev}
p=\frac{qn}{q\eta+n-a},
\end{equation}
we have that the solution $\theta(t,\cdot)$ belongs to a H\"older space $\mathcal{C}^{\gamma}(\mathbb{R}^{n})$ for $0<\gamma<\frac{1}{4}$ small. 
\end{Theoreme}
Just as in the two previous theorems, we have at our disposal two types of information. On the one hand, we have a general Morrey-Campanato boundedness property that reads as follows in this setting:
\begin{equation}
\label{HLS}
\begin{split}
\|A_{[(-\Delta)^{-\frac{\eta}{2}}\theta]}\|_{L^{\infty}(M^{q,a})}&\leq C_{A}\|(-\Delta)^{-\frac{\eta}{2}}\theta\|_{L^{\infty}(M^{q,a})}\leq C_{A}\|(-\Delta)^{-\frac{\eta}{2}}\theta\|_{L^{\infty}(L^{r})}\\
&\leq C_{A}\|\theta\|_{L^{\infty}(L^{p})},
\end{split}
\end{equation}
where we used the Morrey-Campanato-Lebesgue embedding (\ref{Injection_MorreyLebesgue}) setting $r=\frac{qn}{n-a} $ and the Hardy-Littlewood-Sobolev inequality $\|(-\Delta)^{-\frac \eta 2}\theta\|_{L^\infty(L^\sigma)}\le c\|\theta\|_{L^\infty(L^p)}, \sigma=\frac{pn}{n-\eta p} $. Solving $\sigma=r$ yields relationship (\ref{RelationMorreySobolev}). Within this framework, applying the previous general results of \cite{DCHSM}, the corresponding smoothness degree $\alpha_{0}$ is given by the relationship (\ref{RelationMorrey}) and satisfies $1<\alpha_{0}<2$. Note in particular that inequality (\ref{HypoMorreyBorne3}) does not convey any regularity property. \\

On the other hand, we have a Besov control and using the nonlinearity of the drift we can see from (\ref{DegreLevy11}) that we obtain a stronger smoothness degree $\alpha$: indeed the fact that $\alpha<1$ is equivalent to $\frac{n-1}{p}<\eta$, and this is possible if the parameter $p$ is big enough. This situation is to be expected as we introduced the regularizing term $(-\Delta)^{-\frac{\eta}{2}}$, note however that if $\eta\longrightarrow 0$ we recover Theorem \ref{Theoreme1}. Moreover, just as in the previous Theorem \ref{Theoreme1}, if $p\longrightarrow +\infty$ we will loose all the Besov information and we will have $\alpha \longrightarrow 1$.
\textcolor{black}{We emphasize as well that, even if we consider a \textit{regularized} information as argument of the non linear operator $A_{[(-\Delta)^{-\frac{\eta}{2}}\theta]} $, we keep some kind of \textit{homogeneity} property imposing a more stringent integrability condition assuming boundedness in $L^{p}(\dot{B}^{\alpha/p+\eta,p}_{p}) $. }\\

Finally,  in our last theorem we consider a more irregular nonlinear drift in the following way $A_{[(-\Delta)^{\frac{\eta}{2}}\theta]}$ for $0<\eta<1$ and we will investigate the corresponding regularity degree $\alpha$. Such type of singularizations have been considered as well in \cite{CCCGW}.

\begin{Theoreme}[Singular Drift]\label{Theoreme4} Let $\theta_{0}\in L^{p}\cap L^{\infty}(\mathbb{R}^{n})$  be an initial data with $n(n-1)\vee 2n<p<+\infty$, let $0<\eta<\frac{1}{2(n-1)}$ be a fixed (small) parameter and consider the nonlinear equation
\begin{equation}\label{Equation02}
\begin{cases}
\partial_t \theta(t,x)-\nabla\cdot(A_{[(-\Delta)^{\frac{\eta}{2}}\theta]}\,\theta)(t,x)+\mathcal{L}^{\alpha}\theta(t,x)=0,\\[4mm]
\theta(0,x)=\theta_{0}(x). 
\end{cases}
\end{equation}
Assume the following points
\begin{itemize}
\item[1)] $\mathcal{L}^{\alpha}$ is a L\'evy-type operator of the form (\ref{DefinitionOperateur}) satisfying hypothesis (\ref{DefKernelLevy}) with a smoothness degree $\alpha$ such that
\begin{equation}\label{DegreLevy12}
1<\alpha=\frac{n+p(1+\eta)}{p+1}<2.
\end{equation} 
\item[2)] The operator $A_{[(-\Delta)^{\frac{\eta}{2}}\theta]}$ that defines the nonlinear drift in equation (\ref{Equation02}) above 
\begin{itemize}
\item[$(a_{4})$] is bounded in Morrey-Campanato spaces 
\begin{equation}\label{HypoMorreyBorne4}
\|A_{[(-\Delta)^{\frac{\eta}{2}}\theta]}\|_{L^{\infty}(M^{p,a})}\leq C_{A}\|\theta\|_{L^{\infty}(L^{p})}\qquad (0\leq a <n),
\end{equation} 
\item[$(b_{4})$] is bounded in Besov spaces $\|A_{[(-\Delta)^\frac{\eta}2\theta]}\|_{L^{p}(\dot{B}^{\alpha/p-\textcolor{black}{\eta},p}_{p})}\leq D_{A}\|(-\Delta)^\frac \eta 2\theta\|_{L^{p}(\dot{B}^{\alpha/p-\textcolor{black}{\eta},p}_{p})}
$.
\end{itemize}
\end{itemize}
Then if $0<\frac{\alpha}{p}-\eta<1$, as long as we have the following condition between the Morrey-Campanato parameter $a$, the regularity degree $\alpha$ and the Lebesgue index $p$ of the initial data:
\begin{equation}\label{MorreyvsBesov12}
0\leq a< \frac{n+p(1-p\eta)}{p+1}
\end{equation}
the solution $\theta(t,\cdot)$  of the equation (\ref{Equation02})  belongs to a H\"older space $\mathcal{C}^{\gamma}(\mathbb{R}^{n})$ for $0<\gamma<\frac{1}{4}$ small and moreover this value of $\alpha$ satisfies $1<\alpha<\alpha_{0}<2$ and we break down the relationship (\ref{RelationMorrey}).\\

However, if we have 
$$ a\geq  \frac{n+p(1-p\eta)}{p+1},$$ 
then the solution $\theta(t,\cdot)$  of the equation (\ref{Equation02}) still belongs to a H\"older space $\mathcal{C}^{\gamma}(\mathbb{R}^{n})$ for $0<\gamma<\frac{1}{4}$ small, but we have now $1<\alpha_{0}\leq \alpha<2$: the Besov regularity information is useless and we do not break down the relationship (\ref{RelationMorrey}).
\end{Theoreme}
This last theorem is very close to Theorem \ref{Theoreme2}, indeed if $\eta\longrightarrow 0$ we recover the results stated there.\\

We remark that the presence of the operator $(-\Delta)^{\frac{\eta}{2}}$ introduces (as expected) an even more singular behavior and this can be observed if we compare the condition (\ref{MorreyvsBesov1}) with (\ref{MorreyvsBesov12}). In this setting the Besov information is relevant only if the Morrey-Campanato space is quite singular and the threshold for the Morrey-Campanato regularity given by identity $a=\frac{n+p(1-p\eta)}{p+1}$ is indeed lower than the one stated in Theorem \ref{Theoreme2}. \textcolor{black}{Observe as well that, since we singularize the drift, we cannot rely anymore on the Hardy-Littlewood Sobolev inequality as in Theorem \ref{Theoreme3} \textcolor{black}{(see estimates (\ref{HLS}))}. This induces to consider a stronger boundedness condition in $(a_4)$. Similarly, since $A_{[(-\Delta)^\frac{\eta}2\theta]} $ is more singular, we alleviate the integrability constraint in $(b_4) $ considering boundedness in $L^{p}(\dot{B}^{\alpha/p-\eta,p}_{p}) $.}\\

With Theorems \ref{Theoreme1}-\ref{Theoreme4} we show how to exploit, from different points of view and in a nonlinear setting, the additional Besov regularity information given by the general formula (\ref{CTR_BESOV_LP}). As expected, this information allows us to improve previous results but it probably seems that the Besov control is too weak to enhance the situation when Morrey-Campanato convey regular properties. 

\begin{Remarque} Although the hypotheses over the Morrey-Campanato norm differ from Theorem \ref{Theoreme1} to Theorem \ref{Theoreme4}, the Besov control stated in terms of the space $\dot{B}^{\alpha/p,p}_{p}$ is \textcolor{black}{essentially} the \emph{same} in all these theorems. This will allow us to give a general treatment for these theorems. 
\end{Remarque}
\mysection{Existence, maximum principle and other properties.}\label{SecExiUnic}
In this section we state results that are needed in order to perform the computations of the following sections. Some of these properties are well known and rather classical, moreover many of these details were already treated in \cite{DCHSM} but, for the sake of completeness, we will explain here some computations, especially where we need to deal with the nonlinear drift and the different hypotheses used in the theorems stated in the section above.  \\

\begin{itemize}
\item[$(i)$] {\bf Existence}. We start studying existence for weak solutions of equation (\ref{EquationUNA}) with initial data $\theta_0\in L^p\cap L^{\infty}(\mathbb{R}^n)$ where $n(n-1)\vee 2n<p<+\infty$. For this we introduce a regularization and using suitable controls we will pass to the limit to obtain weak solutions. Thus, for  $\varepsilon>0 $, we consider the following problem
\begin{equation}\label{EquationRegular}
\left\lbrace
\begin{array}{l}
\partial_t \theta(t,x)- \nabla\cdot(A^{\varepsilon}_{[\theta]}\;\theta)(t,x)+\mathcal{L}^{\alpha}\theta(t,x)=\varepsilon \Delta \theta(t,x),\qquad t\in [0,T], \\[4mm]
\theta(0,x)=\theta_0(x).
\end{array}
\right.
\end{equation}
To obtain regularity in the time variable, we introduce $A_{[\theta]}^{\star,\varepsilon}=A_{[\theta]}\star \psi_{\varepsilon}$ where $\star$ stands for the time convolution and $\psi_{\varepsilon}(t)=\varepsilon^{-1}\psi(t/\varepsilon)$ where  $\psi \in \mathcal{C}^{\infty}_0(\mathbb{R})$ is a non-negative function such that $supp(\psi)\subset B(0,1)$ and $\displaystyle{\int_{\mathbb{R}}}\psi(t)dt=1$. In the previous time convolution, we have extended the velocity field (in the time variable) on $\R $ setting for all $s\in \R\setminus [0,T]$, $A_{[\theta]}(s,\cdot)=0 $.
Then we define $A_{[\theta]}^{\varepsilon}=A_{[\theta]}^{\star,\varepsilon}\ast \omega_{\varepsilon}$,  here  $*$ stands now for the spatial convolution and $\omega_{\varepsilon}$ is a usual mollifying kernel, i.e. $\omega_{\varepsilon}(x)=\varepsilon^{-n}\omega(x/\varepsilon)$, $\omega \in \mathcal{C}^{\infty}_0(\mathbb{R}^n)$ is a non-negative function such that $supp(\omega)\subset B(0,1)$ and $\displaystyle{\int_{\mathbb{R}^n}}\omega(x)dx=1$. From this regularization, for a fixed $\varepsilon>0$, the approximate drift $A_{[\theta]}^{\varepsilon}$ will be a smooth (in the time and space variables), divergence free vector field. \\

Since we can apply this approximation procedure in a completely similar manner to the drifts $A_{[(-\Delta)^{-\frac{\eta}{2}}\theta]}$ and $A_{[(-\Delta)^{\frac{\eta}{2}}\theta]}$ used in Theorems \ref{Theoreme3} and \ref{Theoreme4} respectively, we will continue our study by considering the problem (\ref{EquationRegular}) which corresponds to the setting of Theorems \ref{Theoreme1} and \ref{Theoreme2}.\\

The problem (\ref{EquationRegular}) can thus be rewritten in the following way:
$$\theta(t,x)=e^{\varepsilon t\Delta}\theta_0(x)+\int_{0}^t e^{\varepsilon (t-s)\Delta}\nabla \cdot(A_{[\theta]}^{\varepsilon}\; \theta)(s,x)ds-\int_{0}^t e^{\varepsilon (t-s)\Delta}\mathcal{L}^{\alpha} \theta(s,x)ds,$$
and, in order to construct solutions, we will apply the Banach contraction scheme in the space $L^{\infty}([0,T], L^{p}(\mathbb{R}^n))$ with the norm $\|\theta\|_{L^\infty (L^p)}=\displaystyle{\underset{t\in [0,T]}{\sup}}\|\theta(t,\cdot)\|_{L^p}$. Following \cite{DCHSM} we have the estimates
\begin{equation}\label{PointFixe1}
\begin{split}
\|e^{\varepsilon t \Delta}\theta_{0}\|_{L^\infty (L^p)}& \leq \|\theta_{0}\|_{L^p} \quad \mbox{and} \\
\left\|\int_{0}^t e^{\varepsilon (t-s)\Delta}\mathcal{L}^{\alpha} \theta(s,x)ds\right\|_{L^\infty (L^p)}&\leq C \left(\frac{T^{1-\frac{\alpha}{2}}}{\varepsilon^{\frac{\alpha}{2}}}+\frac{T^{1-\frac{\delta}{2}}}{\varepsilon^{\frac{\delta}{2}}}\right)\; \|\theta\|_{L^\infty (L^p)},
\end{split}
\end{equation}
where we recall from  \eqref{DefKernelLevy} that the exponent $\delta $ is related to the tail behavior of the operator ${\mathcal L}^\alpha $.

Thus, we only need to treat the bilinear quantity
$$B(\theta, \varphi)=\int_{0}^t e^{\varepsilon (t-s)\Delta}\nabla \cdot(A_{[\theta]}^{\varepsilon}\; \varphi)(s,x)ds,$$
and we have
\begin{eqnarray}
\|B(\theta, \varphi)\|_{L^\infty (L^p)}& =&\underset{0<t<T'}{\sup} \left\|\int_{0}^t e^{\varepsilon (t-s)\Delta} \nabla \cdot(A_{[\theta]}^{\varepsilon}\; \varphi)(s,\cdot)ds\right\|_{L^p}\nonumber\\
&=&\underset{0<t<T'}{\sup} \left\|\int_{0}^t \nabla \cdot(A_{[\theta]}^{\varepsilon}\; \varphi)\ast h_{\varepsilon (t-s)}(s,\cdot)ds\right\|_{L^p}\nonumber\\
&\leq & \underset{0<t<T'}{\sup}\int_{0}^t  \left\|(A_{[\theta]}^{\varepsilon}\; \varphi)(s,\cdot)\right\|_{L^p} \left\|\nabla h_{\varepsilon (t-s)}\right\|_{L^1} ds\nonumber\\
&\leq & \underset{0<t<T'}{\sup}\int_{0}^t  \left\|A_{[\theta]}^{\varepsilon}(s,\cdot)\right\|_{L^\infty}  \left\|\varphi(s,\cdot)\right\|_{L^p} C(\varepsilon(t-s))^{-1/2} ds.\label{PointFixe11}
\end{eqnarray}
The computations performed until now (\emph{i.e.} from (\ref{EquationRegular}) to (\ref{PointFixe11})) are common to Theorems \ref{Theoreme1}-\ref{Theoreme4} and we will continue the proof of the existence under the hypotheses of Theorem \ref{Theoreme1}: the modifications for Theorems \ref{Theoreme2}, \ref{Theoreme3} and \ref{Theoreme4} are given in Remark \ref{RemarqueExistence} below.\\

At this point we use the following inequality 
\begin{equation}\label{EstimationMorreyUtile}
\left\|A_{[\theta]}^{\varepsilon}
\right\|_{L^\infty(L^\infty)}\leq C\varepsilon^{-n/q} \|A_{[\theta]}\|_{L^\infty (M^{q,a})},
\end{equation}
 valid for $1<q<+\infty$ and which is proven in the Lemma A-1 of \cite{DCHSM} \textcolor{black}{(for functions in  $M^{q,a} $)} to obtain
\begin{eqnarray*}
\|B(\theta, \varphi)\|_{L^\infty (L^p)}&\leq & C\varepsilon^{-n/q} \|A_{[\theta]}\|_{L^\infty (M^{q,a})}\left\|\varphi\right\|_{L^\infty (L^p)}  \underset{0<t<T'}{\sup}\int_{0}^t C(\varepsilon(t-s))^{-1/2} ds\nonumber\\
&\leq & C \frac{T^{\frac 12}}{\varepsilon^{\frac 12}}\;  \varepsilon^{-n/q}\|A_{[\theta]}\|_{L^\infty (M^{q,a})}  \left\|\varphi\right\|_{L^\infty (L^p)}.
\end{eqnarray*}
We use now hypothesis (\ref{HypoMorreyBorne1}) stated in Theorem \ref{Theoreme1} (\emph{i.e.} $\|A_{[\theta]}\|_{L^\infty (M^{q,a})} \leq C_{A}\|\theta\|_{L^\infty (M^{q,a})}$) and we have
\begin{equation}\label{PointFixe2}
\|B(\theta, \varphi)\|_{L^\infty (L^p)}\leq C_{A} \frac{T^{\frac 12}}{\varepsilon^{\frac 12}}\;\varepsilon^{- n/q}\|\theta\|_{L^\infty (M^{q,a})} \|\varphi\|_{L^\infty (L^p)},
\end{equation}
but, since by the relationship (\ref{Equation_Existence}) we have the identity
$p=\frac{qn}{n-a}$, we use the embedding (\ref{Injection_MorreyLebesgue}) to write $\|\theta\|_{L^\infty (M^{q,a})}\leq \|\theta\|_{L^\infty (L^{p})}$, from which we deduce the continuity of the bilinear form $B(\theta, \varphi)$:
\begin{equation}\label{Formula_FormaBilineal}
\|B(\theta, \varphi)\|_{L^\infty (L^p)}\leq C_{A} \frac{T^{\frac 12}}{\varepsilon^{\frac 12}}\;\varepsilon^{- n/q}\|\theta\|_{L^\infty (L^{p})} \|\varphi\|_{L^\infty (L^p)}.
\end{equation}
 
Now, with inequalities (\ref{PointFixe1}) and (\ref{Formula_FormaBilineal}), it is easy to obtain the existence of solutions for the regularized problem (\ref{EquationRegular}) in a small time interval.
\begin{Remarque}[Existence for Theorems \ref{Theoreme2}, \ref{Theoreme3} \& \ref{Theoreme4}]\label{RemarqueExistence}
\begin{itemize}
\item[]
\item[$\bullet$] To prove the existence of solutions for the problem (\ref{EquationRegular}) under the hypotheses stated in Theorem \ref{Theoreme2}, we only need to study inequality (\ref{PointFixe11}) and for this we use the general estimate (\ref{EstimationMorreyUtile}) in the following manner $\left\|A_{[\theta]}^{\varepsilon}(s,\cdot)\right\|_{L^\infty(L^\infty)}\leq C\varepsilon^{-n/p} \|A_{[\theta]}\|_{L^\infty (M^{p,a})}$, 
thus, applying the hypothesis (\ref{HypoMorreyBorne2}) we obtain the following inequality
$$\|B(\theta, \varphi)\|_{L^\infty (L^p)}\leq C_{A} \frac{T^{\frac 12}}{\varepsilon^{\frac 12}}\;\varepsilon^{- n/p}\|\theta\|_{L^\infty (L^{p})} \|\varphi\|_{L^\infty (L^p)},$$
from which we deduce the continuity of the bilinear form. 
\item[$\bullet$] For the existence of solutions for equation (\ref{EquationRegular}) within the framework of Theorem \ref{Theoreme3} we start from expression (\ref{PointFixe11}) and we have
$$\|B(\theta, \varphi)\|_{L^\infty (L^p)}\leq \underset{0<t<T'}{\sup}\int_{0}^t  \left\|A_{[(-\Delta)^{-\frac{\eta}{2}}\theta]}^{\varepsilon}(s,\cdot)\right\|_{L^\infty}  \left\|\varphi(s,\cdot)\right\|_{L^p} C(\varepsilon(t-s))^{-1/2} ds,$$
thus using inequality (\ref{EstimationMorreyUtile}) we obtain the continuity of this bilinear form:
\begin{eqnarray*}
\|B(\theta, \varphi)\|_{L^\infty (L^p)}&\leq & C \frac{T^{\frac 12}}{\varepsilon^{\frac 12}}\;  \varepsilon^{-n/q}\|A_{[(-\Delta)^{-\frac{\eta}{2}}\theta]}\|_{L^\infty (M^{q,a})}  \left\|\varphi\right\|_{L^\infty (L^p)}\\
&\leq & C_{A} \frac{T^{\frac 12}}{\varepsilon^{\frac 12}}\;  \varepsilon^{-n/q}\|\theta\|_{L^\infty (L^{p})}  \left\|\varphi\right\|_{L^\infty (L^p)},
\end{eqnarray*}
where we used \eqref{HLS} for the last inequality.

\item[$\bullet$] Finally, for Theorem \ref{Theoreme4}, from expression (\ref{PointFixe11}) we have
$$\|B(\theta, \varphi)\|_{L^\infty (L^p)}\leq \underset{0<t<T'}{\sup}\int_{0}^t  \left\|A_{[(-\Delta)^{\frac{\eta}{2}}\theta]}^{\varepsilon}(s,\cdot)\right\|_{L^\infty}  \left\|\varphi(s,\cdot)\right\|_{L^p} C(\varepsilon(t-s))^{-1/2} ds,$$ applying the general inequality (\ref{EstimationMorreyUtile}) we obtain
$$\left\|A_{[(-\Delta)^{\frac{\eta}{2}}\theta]}^{\varepsilon}(s,\cdot)\right\|_{L^\infty(L^\infty)}\leq C\varepsilon^{-n/p} \|A_{[(-\Delta)^{\frac{\eta}{2}}\theta]}\|_{L^\infty (M^{p,a})},$$ 
and with the hypothesis (\ref{HypoMorreyBorne4}) we write
\begin{eqnarray*}
\|B(\theta, \varphi)\|_{L^\infty (L^p)}&\leq & C \frac{T^{\frac 12}}{\varepsilon^{\frac 12}}\;  \varepsilon^{-n/p}\|A_{[(-\Delta)^{\frac{\eta}{2}}\theta]}\|_{L^\infty (M^{p,a})}  \left\|\varphi\right\|_{L^\infty (L^p)}\\
&\leq & C_{A} \frac{T^{\frac 12}}{\varepsilon^{\frac 12}}\;  \varepsilon^{-n/p}\|\theta\|_{L^\infty (L^{p})}  \left\|\varphi\right\|_{L^\infty (L^p)}.
\end{eqnarray*}
\end{itemize}
\end{Remarque}

\item[$(ii)$] {\bf Regularity}. We now study the regularity of the solutions constructed by this method. And we will see that solutions of the approximated problem (\ref{EquationRegular}) are actually smooth. Indeed, following \cite{DCHSM}, we will work in the time interval $0<T_{0}<T_{\ast}<t<T^{\ast}$ where $T_{0}$, $T_{\ast}$ and $T^{\ast}$ are fixed bounds and by iteration we will prove that $\theta \in \displaystyle{\bigcap_{0<T_0<T_{\ast}<t<T^\ast}}L^\infty([0,t], W^{\frac{k}{2},p}(\mathbb{R}^n))$ for all $k\in {\mathbb N}$ where we define the Sobolev space $W^{s,p}(\mathbb{R}^n)$ for $ s\in \mathbb{R}$ and $1<p<+\infty$ by  $\|f\|_{W^{s,p}}=\|f\|_{L^p}+\|(-\Delta)^{\frac{s}{2}}f\|_{L^p}$. Note that this is true for $k=0$. So let us assume that it is also true for $k>0$ and we will show that it is still true for $k+1$.\\

We will work under the framework given in Theorem \ref{Theoreme1}, modifications for Theorems \ref{Theoreme2}-\ref{Theoreme4} are straightforward using the computations made in Remark \ref{RemarqueExistence}.\\

We thus consider then the next problem
$$\theta(t,x)=e^{\varepsilon (t-T_0)\Delta}\theta(T_0,x)+\int_{T_0}^t e^{\varepsilon (t-s)\Delta}\nabla \cdot(A^{\varepsilon}_{[\theta]}\; \theta)(s,x)ds-\int_{T_0}^t e^{\varepsilon (t-s)\Delta}\mathcal{L}^{\alpha} \theta(s,x)ds,
$$
and we have the following estimate
\begin{equation}\label{InegaliteRegulariteApprochee}
\begin{split}
\|\theta\|_{L^\infty (W^{\frac{k+1}{2},p})}&\leq \|e^{\varepsilon (t-T_0)\Delta}\theta(T_0,\cdot)\|_{L^\infty (W^{\frac{k+1}{2},p})}\\
& +\left\|\int_{T_0}^t e^{\varepsilon (t-s)\Delta}\nabla \cdot(A^{\varepsilon}_{[\theta]}\; \theta)(s,\cdot)ds\right\|_{L^\infty (W^{\frac{k+1}{2},p})}\\
& +\left\|\int_{T_0}^t e^{\varepsilon (t-s)\Delta}\mathcal{L}^{\alpha}\theta(s,\cdot)ds\right\|_{L^\infty (W^{\frac{k+1}{2},p})}.
\end{split}
\end{equation}
The first and the last term of (\ref{InegaliteRegulariteApprochee}) were already treated in \cite{DCHSM}, indeed we have
\begin{eqnarray*}
\|e^{\varepsilon (t-T_0)\Delta}\theta(T_0,\cdot)\|_{L^\infty (W^{\frac{k+1}{2},p})}&\leq &C\|\theta(T_0,\cdot)\|_{L^p}\max\left\{1; \left[\varepsilon (t-T_0)\right]^{- \frac{k+1}{4}}\right\}\qquad \mbox{ and }\\
\left\|\int_{T_0}^t e^{\varepsilon (t-s)\Delta} \mathcal{L}  \theta(s,\cdot)ds\right\|_{L^\infty (W^{\frac{k+1}{2},p})}&\leq & C\|\theta \|_{L^{\infty}(W^{\frac{k}{2},p})}
\end{eqnarray*}
$$\times \int_{T_0}^t \max\left\{\left([\varepsilon(t-s)]^{-\frac{\alpha}{2}}+\varepsilon(t-s)]^{-\frac{\delta}{2}}\right) ; \big([\varepsilon (t-s)]^{-\frac{1+2\alpha}{4}}+[\varepsilon (t-s)]^{-\frac{1+2\delta}{4}} \big)   \right\}ds.$$
Thus we only study the second term of (\ref{InegaliteRegulariteApprochee}) and following the same ideas as in the existence computations performed above one has
\begin{equation}\label{RegulariteTheo}
\begin{split}
I=\left\|\int_{T_0}^t e^{\varepsilon (t-s)\Delta}\nabla \cdot(A^{\varepsilon}_{[\theta]}\; \theta)(s,\cdot)ds\right\|_{L^\infty (W^{\frac{k+1}{2},p})}&\leq C \|A_{[\theta]}\|_{L^\infty (M^{q,a})}\|\theta\|_{L^\infty (W^{\frac{k}{2},p})}\\
 \times \sup \int_{T_0}^t  \left(1+\varepsilon^{-N}\right)\varepsilon^{-n/q}&\max\left\{  \left[\varepsilon (t-s)\right]^{- \frac{1}{2}} ;  \left[\varepsilon (t-s)\right]^{- \frac{3}{4}}\right\}ds,
\end{split}
\end{equation}
for $N\ge \frac k2 $.
At this point we use the hypotheses of Theorem \ref{Theoreme1} (\emph{i.e.} the Morrey-Campanato boundedness and the Morrey-Campanato-Lebesgue embedding) to write
\begin{equation}\label{RegulariteTheo1}
\begin{split}
I &\leq C_{A} \|\theta\|_{L^\infty  (M^{q,a})}\|\theta\|_{L^\infty (W^{\frac{k}{2},p})} \\
& \times\sup \int_{T_0}^t  \left(1+\varepsilon^{-N}\right)\varepsilon^{-n/q}\max\left\{  \left[\varepsilon (t-s)\right]^{- \frac{1}{2}} ;  \left[\varepsilon (t-s)\right]^{- \frac{3}{4}}\right\}ds\\
&\leq C_{A} \|\theta\|_{L^\infty  (L^{p})}\|\theta\|_{L^\infty (W^{\frac{k}{2},p})} \\
& \times \sup \int_{T_0}^t  \left(1+\varepsilon^{-N}\right)\varepsilon^{-n/q}\max\left\{  \left[\varepsilon (t-s)\right]^{- \frac{1}{2}} ;  \left[\varepsilon (t-s)\right]^{- \frac{3}{4}}\right\}ds.
\end{split}
\end{equation}

With these inequalities we conclude that the norm $\|\theta\|_{L^\infty (W^{\frac{k+1}{2},p})}$ is controlled for all $\varepsilon>0$: we have proven spatial regularity. Time regularity inductively follows since we have for all $ k\geq 0$:
$$\frac{ \partial^{k+1}}{\partial t^{k+1}}\theta(t,x)-\nabla \cdot \left(\frac{ \partial^k}{\partial t^k} (A^{\varepsilon}_{[\theta]}\,\theta)\right)(t,x)+\mathcal{L}^{\alpha}\left(\frac{ \partial^k}{\partial t^k} \theta\right)(t,x)=\varepsilon \Delta \left(\frac{ \partial^k}{\partial t^k} \theta\right)(t,x).$$

Regularity under the hypotheses of Theorems \ref{Theoreme2}-\ref{Theoreme4} follows essentially the same lines: apply Remark \ref{RemarqueExistence}
to the computations performed in (\ref{RegulariteTheo})-(\ref{RegulariteTheo1}).

 \item[$(iii)$] {\bf Maximum principle}.\label{MaxPrinciple} Let $\theta_0\in L^{p}(\mathbb{R}^n)$ with $n(n-1)\vee 2n< p< +\infty$ be an initial data, then the associated solution of the regularized problem  (\ref{Equation02}) satisfies the following maximum principle for all $t\in [0, T]$: 
\begin{equation}\label{MaximumPrinciple}
\|\theta(t,\cdot)\|_{L^p}\leq \|\theta_0\|_{L^p}.
\end{equation}
This maximum principle is a consequence of the divergence free property of the drift and of the following computations. Again, we assume the framework of Theorem \ref{Theoreme1} since the result for Theorems \ref{Theoreme2}-\ref{Theoreme4} is essentially the same.  We write for $n(n-1)\vee 2n<p<+\infty$:
\begin{eqnarray*}
\frac{d}{dt}\|\theta(t,\cdot)\|^p_{L^p}&=&p\int_{\mathbb{R}^n}|\theta|^{p-2}\theta\bigg(\varepsilon \Delta \theta+\nabla \cdot (A^{\varepsilon}_{[\theta]} \,\theta)-\mathcal{L}^{\alpha}\theta\bigg)(t,x)dx\\
&=&p\varepsilon\int_{\mathbb{R}^n}|\theta|^{p-2}\theta\Delta \theta(t,x) dx-p\int_{\mathbb{R}^n}|\theta|^{p-1}sgn(\theta) \mathcal{L}^{\alpha}\theta(t,x) dx,
\end{eqnarray*}
where we used the fact that $\nabla \cdot (A^{\varepsilon}_{[\theta]})=0$. Thus, we have
$$\frac{d}{dt}\|\theta(t,\cdot)\|^p_{L^p}-p\varepsilon\int_{\mathbb{R}^n}|\theta|^{p-2}\theta \Delta \theta(t,x)dx+p\int_{\mathbb{R}^n}|\theta|^{p-1}sgn(\theta)\mathcal{L}^{\alpha}\theta (t,x) dx=0,$$
and integrating in time we obtain
\begin{equation}\label{PrincipeduMax1}
\begin{split}
\|\theta(t,\cdot)\|^p_{L^p}&-p\varepsilon\int_{0}^t\int_{\mathbb{R}^n}|\theta|^{p-2}\theta\Delta \theta(s,x) dxds\\
&+p\int_0^t\int_{\mathbb{R}^n}|\theta|^{p-1}sgn(\theta)\mathcal{L}^{\alpha} \theta (s,x) dxds=\|\theta_0\|^p_{L^p}.
\end{split}
\end{equation}
To finish, we have that the quantities 
\begin{equation*}
-p\varepsilon\displaystyle{\int_{\mathbb{R}^n}}|\theta|^{p-2}\theta \Delta \theta(s,x) dx\quad \mbox{ and }\quad \displaystyle{p\int_0^t\int_{\mathbb{R}^n}}|\theta|^{p-1}sgn(\theta) \mathcal{L}^{\alpha} \theta(s,x) dxds
\end{equation*}
are both positive. We refer to the prof of Proposition 2.1 in \cite{DCHSM} for details.  Thus, getting back to (\ref{PrincipeduMax1}), we have that all these quantities are bounded and positive and we obtain the maximum principle property stated in (\ref{MaximumPrinciple}).
\begin{Remarque}
With this a priori bound (\ref{MaximumPrinciple}), it is posible to obtain global solutions $\theta(t, \cdot)$ where $t\in [0,T]$ and $0<T<+\infty$ is a fixed final time. 
\end{Remarque}

\begin{Remarque}[Uniform controls]\label{RemarqueControlUniformeTransport} Following the hypotheses stated in Theorems \ref{Theoreme1}-\ref{Theoreme4}, this maximum principle implies a uniform control over the nonlinear drift. 
\begin{itemize}
\item[$\bullet$] for Theorem \ref{Theoreme1}, using (\ref{HypoMorreyBorne1}) and the Morrey-Campanato-Lebesgue embedding (\ref{Injection_MorreyLebesgue}) we have
$$\|A_{[\theta]}\|_{L^{\infty}(M^{q,a})}\leq C_{A}\|\theta\|_{L^{\infty}(M^{q,a})}\leq C_{A}\|\theta\|_{L^{\infty}(L^{p})}\leq C_{A}\|\theta_0\|_{L^p}<+\infty,$$
\item[$\bullet$] for Theorem \ref{Theoreme2}, using (\ref{HypoMorreyBorne2}) we have
$$\|A_{[\theta]}\|_{L^{\infty}(M^{p,a})}\leq C_{A}\|\theta\|_{L^{\infty}(L^{p})}\leq C_{A}\|\theta_0\|_{L^p}<+\infty,$$
\item[$\bullet$] for Theorem \ref{Theoreme3}, from \eqref{HLS},
\begin{eqnarray*}
\|A_{[(-\Delta)^{-\frac{\eta}{2}}\theta]}\|_{L^{\infty}(M^{q,a})}
\leq C_{A}\|\theta\|_{L^{\infty}(L^{p})} \leq  C_{A}\|\theta_0\|_{L^p}<+\infty,
\end{eqnarray*}
\item[$\bullet$] finally, for Theorem \ref{Theoreme4}, using (\ref{HypoMorreyBorne4}) we have
$$\|A_{[(-\Delta)^{\frac{\eta}{2}}\theta]}\|_{L^{\infty}(M^{p,a})}\leq C_{A}\|\theta\|_{L^{\infty}(L^{p})}\leq C_{A}\|\theta_0\|_{L^p}<+\infty.$$

\end{itemize}
\end{Remarque}
\item[$(iv)$] {\bf Besov Regularity}. This property is the crucial point of this article as it gives an additional information on weak solutions which was not exploited in \cite{DCHSM} where we only considered linear equations. Thus, since we consider now a nonlinear drift that satisfies a Besov space stability property as stated in Theorems \ref{Theoreme1}-\ref{Theoreme4}, it is interesting to investigate the behavior of weak solutions in terms of Besov spaces and to this end we have the following theorem from \cite{PGDCH}:
\begin{Theoreme}[Besov Regularity]\label{RegulariteBesov} Let $\mathcal{L}^{\alpha}$ be a Lévy-type operator of the form (\ref{DefinitionOperateur}) satisfying hypothesis (\ref{DefKernelLevy}) with $0<\alpha<2$.  Let $2\leq p <+\infty$ and let $f:\mathbb{R}^n\longrightarrow \mathbb{R}$ be a function such that $f\in L^{p}(\mathbb {R}^n)$ and
$$\int_{\mathbb{R}^n}|f(x)|^{p-2}f(x)\mathcal{L}^{\alpha} f(x)dx<+\infty, \quad\mbox{ then  }\quad f\in \dot{B}^{\frac{\alpha}{p},p}_{p}(\mathbb{R}^n).$$
\end{Theoreme}
See \cite{PGDCH} and \cite{DCHSM} for a proof of this fact. Let us explain briefly how we will use this regularity property: \textcolor{black}{combining this result and formula (\ref{PrincipeduMax1}) we obtain the important estimate \eqref{CTR_BESOV_LP} which gives $\theta\in L^p([0,T],\dot{B}^{\frac{\alpha}{p},p}_{p}) =L^p(\dot{B}^{\frac{\alpha}{p},p}_{p})$ and therefore that for almost all $t\in [0,T] $, $\theta(t,\cdot) \in \dot{B}^{\frac{\alpha}{p},p}_{p}$.} 
Now, from the Besov stability hypothesis stated in Theorems \ref{Theoreme1}-\ref{Theoreme4}  we have $\|A_{[\theta]}\|_{L^{p}(\dot{B}^{\frac{\alpha}{p},p}_{p})}\leq D_{A}\|\theta\|_{L^{p}(\dot{B}^{\frac{\alpha}{p},p}_{p})}$, 
which gives us a Besov control for the nonlinear drift:
\begin{equation}\label{FormuleControlMorreyBesov}
\|A_{[\theta]}\|_{L^{p}(\dot{B}^{\frac{\alpha}{p},p}_{p})}\leq D_{A}\|\theta\|_{L^{p}(\dot{B}^{\frac{\alpha}{p},p}_{p})}\leq D_{A}\|\theta_{0}\|_{L^{p}}<+\infty.
\end{equation}

\item[$(v)$] {\bf The limit}. The previous solutions form a family of regular functions $(\theta^{(\varepsilon)})_{\varepsilon >0}\in L^{\infty}([0,T], L^{p}(\mathbb{R}^n))$ which satisfy the uniform bound $\|\theta^{(\varepsilon)}(t,\cdot)\|_{L^p}\leq \|\theta_0\|_{L^p}$; in order to conclude we need to pass to the limit by letting $\varepsilon \longrightarrow 0$. For this, since we have $\theta_{0}\in L^{p}\cap L^{\infty}(\mathbb{R}^n)$ and since we have the uniform controls over the nonlinear drift given in Remark \ref{RemarqueControlUniformeTransport}, we can apply the same arguments used in \cite{DCHSM}, we thus obtain existence and uniqueness of weak solutions for the equations stated in Theorems \ref{Theoreme1}-\ref{Theoreme4} with an initial data in $\theta_0\in L^p(\mathbb{R}^n)$, $2n\leq p<+\infty$ that satisfy the maximum principle. Moreover, we have that these solutions $\theta(t,x)$ belong to the space $L^{\infty}([0,T], L^{p}(\mathbb{R}^n))\cap L^{p}([0,T], \dot{B}^{\frac{\alpha}{p},p}_p(\mathbb{R}^n))$.
\begin{Remarque} The hypothesis $\theta_{0}\in L^{\infty}(\mathbb{R}^n)$ (and not only  $\theta_{0}\in L^{p}$) is useful in order to obtain general existence results when dealing with general drifts belonging to Morrey-Campanato spaces. See Theorem 1 of \cite{DCHSM} for the details. 
\end{Remarque}
\end{itemize}


\mysection{H\"older Regularity and proof for Theorems \ref{Theoreme1}-\ref{Theoreme4}}\label{SeccHolderRegularity}
One important feature of Hardy spaces relies in the following property:  let $\frac{n}{n+1}<\sigma<1$ and fix $0<\gamma=n(\frac{1}{\sigma}-1)<1$, then the dual of the local Hardy space $h^{\sigma}(\mathbb{R}^n)$ is the Hölder space\footnote{See \cite{Gold} for a proof of this fact and see \cite{Coifmann} and \cite{Stein2} for a detailed treatment on Hardy spaces.} $\mathcal{C}^\gamma(\mathbb{R}^n)$ and with the help of a dual equation we will use this property in order to study H\"older regularity of the solutions of equation (\ref{EquationUNA}).\\

We recall that the local Hardy space $h^{\sigma}(\mathbb{R}^n)$ with $0<\sigma<1$ is the set of distributions $f$ that admits the decomposition $f=\displaystyle{\sum_{j\in \mathbb{N}}}\lambda_j \psi_j$, where $(\lambda_j)_{j\in \mathbb{N}}$ is a sequence of complex numbers such that $\displaystyle{\sum_{j\in \mathbb{N}}}|\lambda_j|^\sigma<+\infty$ and $(\psi_j)_{j\in \mathbb{N}}$ is a family of molecules given by the following definition.
\begin{Definition}\label{DefMolecules} Set $\frac{n}{n+1}<\sigma<1$, define $\gamma=n(\frac{1}{\sigma}-1)$ and fix a real number $\omega$ such that $0<\gamma<\omega<\frac{1}{4}$. For $0<r<1$ an integrable function $\psi$ is an  molecule if we have
\begin{eqnarray}
& & \int_{\mathbb{R}^n} |\psi(x)||x-x_0|^{\omega}dx \leq  (\zeta r)^{\omega-\gamma}\mbox{, for } x_0\in \mathbb{R}^n, \qquad \mbox{and }\qquad \|\psi\|_{L^\infty}  \leq  \frac{1}{(\zeta r)^{n+\gamma}}\label{Hipo2},\\[1mm]
& &\int_{\mathbb{R}^n} \psi(x)dx=0\label{Hipo3}.
\end{eqnarray}
In the above conditions the quantity $\zeta>1$ denotes a constant that depends on $\gamma, \omega, \alpha$ and other parameters to be specified later on. In the case when $1\leq r <+\infty$, we only require conditions (\ref{Hipo2}) for the  molecule $\psi$ while the moment condition (\ref{Hipo3}) is dropped.
\end{Definition}
\begin{Remarque}\label{Remark2}
Conditions (\ref{Hipo2}) imply the estimate $\|\psi\|_{L^1}\leq C\, (\zeta r)^{-\gamma}$. Thus, every molecule belongs to $L^p(\mathbb{R}^n)$ with $1<p<+\infty$ since we have
\begin{equation}\label{Hipo4}
\|\psi\|_{L^p} \leq C (\zeta r)^{-n+\frac{n}{p}-\gamma}.
\end{equation}
\textcolor{black}{Recall also that the Schwartz class $\mathcal{S}(\mathbb{R}^{n})$ is dense in $h^{\sigma}(\mathbb{R}^{n})$, this fact is of course very useful for approximation procedures.}
\end{Remarque}
See \cite{DCHSM} for more remarks about this definition. See also \cite{Stein2}, Chapter III,  \S 5.7,  \cite{Torchinski}, Chapter XIV, \S 6.6 or \cite{KN}.
\subsection*{The dual equation}
Once we have described the elements of the Hardy spaces, we need to derive a dual equation from the original problem (\ref{EquationUNA}), and for this we proceed as follows: assume that $\theta_{0}\in L^{p}\cap L^{\infty}(\mathbb{R}^{n})$ is an initial data and consider the equation (\ref{EquationUNA}) with a L\'evy-type operator $\mathcal{L}^{\alpha}$ and a nonlinear drift $A_{[\theta]}$ satisfying all the hypotheses stated in Theorem \ref{Theoreme1} (since all the nonlinear drift in Theorems \ref{Theoreme2}-\ref{Theoreme4} are divergence free, the arguments are similar in these cases). By the previous sections we can construct on the interval $[0,T]$, with $0<T<+\infty$ fixed, a corresponding weak solution $\theta(\cdot,\cdot)\in L^{\infty}\big([0,T], L^{p}(\mathbb{R}^{n})\big)\cap L^{p}([0,T], \dot{B}^{\frac{\alpha}{p},p}_p(\mathbb{R}^n))$ that satisfies the maximum principle from which we deduce the uniform bounds (recall Remark \ref{RemarqueControlUniformeTransport} and formula (\ref{FormuleControlMorreyBesov})):
\begin{equation}\label{ControlUniformes}
\begin{split}
\|A_{[\theta]}\|_{L^{\infty}(M^{q,a})}&<+\infty,\\[3mm]
\|A_{[\theta]}\|_{L^{p}(\dot{B}^{\frac{\alpha}{p},p}_p)}&\leq D_{A}\|\theta\|_{L^{p}(\dot{B}^{\frac{\alpha}{p},p}_p)}\leq D_{A}\| \theta_{0}\|_{L^{p}}<+\infty.
\end{split}
\end{equation}
Fix now a time $0<t<T$ and consider the following dual equation
\begin{equation}\label{Evolution01}
\begin{cases}
\partial_s \psi(s,x)=& -\nabla\cdot [A_{[\theta]}(t-s,x)\psi(s,x)]-\mathcal{L}^{\alpha}\psi(s,x),\qquad s\in [0,t],\\[2mm]
\psi(0,x)=& \psi_0(x)\in L^{1}\cap L^{\infty}(\mathbb{R}^n),
\end{cases}
\end{equation}
where $\mathcal{L}^{\alpha}$ is a L\'evy-type operator of degree $0<\alpha<2$ and the initial data $\psi_{0}$ is a  molecule in the sense of Definition \ref{DefMolecules} above.\\ 

Note that this problem is now a \emph{linear} equation since the drift $A_{[\theta]}$ does not depend on $\psi(t,x)$. Remark also that, by hypothesis we have $div(A_{[\theta]})=0$ and from (\ref{ControlUniformes}) we have a uniform control over $\|A_{[\theta]}\|_{L^{\infty}(M^{q,a})}$ and thus by \cite{DCHSM} we obtain global existence and a $L^{p}$-maximum principle for the solutions $\psi(\cdot, \cdot)$ of such equation.\\ 

Once we have the dual equation and the corresponding solution $\psi$, we can state the following theorem which will help us to prove Theorems \ref{Theoreme1}-\ref{Theoreme4} using the Hardy-H\"older duality.
\begin{Theoreme}[Transfer Property]\label{Transfert} Let $\theta_{0}\in L^{p}\cap L^{\infty}(\mathbb{R}^{n})$ and let $\theta(t,x)$ be a weak solution of the equation (\ref{EquationUNA}) in the interval $[0,T]$ where the nonlinear drift and the L\'evy-type operator $\mathcal{L}^{\alpha}$ satisfy the hypotheses of Theorems \ref{Theoreme1}-\ref{Theoreme4}. Let $t\in [0,T]$ and let $\psi(s, x)$ be a solution of the backward problem (\ref{Evolution01}) for $s\in [0,t]$. Then we have the identity
\begin{equation*}
\int_{\mathbb{R}^n}\theta(t,x)\psi(0,x)dx=\int_{\mathbb{R}^n}\theta\left(\frac{t}{2},x\right)\psi\left(\frac{t}{2},x\right)dx.
\end{equation*}
\end{Theoreme}
\textit{\textbf{Proof.}}
We first consider the expression
$$\partial_s\int_{\mathbb{R}^n}\theta(t-s,x)\psi(s,x)dx=\int_{\mathbb{R}^n}-\partial_t\theta(t-s,x)\psi(s,x)+\partial_s\psi(s,x)\theta(t-s,x)dx.$$
Using equations (\ref{EquationUNA}) and (\ref{Evolution01}) we obtain
\begin{eqnarray*}
\partial_s\int_{\mathbb{R}^n}\theta(t-s,x)\psi(s,x)dx&=&\int_{\mathbb{R}^n}\bigg\{ -\nabla\cdot\left[(A_{[\theta]}(t-s,x)\theta(t-s,x)\right]\psi(s,x)+\mathcal{L}^{\alpha}\theta(t-s,x)\psi(s,x) \\[5mm]
&& -\nabla\cdot\left[(A_{[\theta]}(t-s,x) \psi(s,x))\right]\theta(t-s,x)-\mathcal{L}^{\alpha}\psi(s,x) \theta(t-s,x)\bigg\}dx.
\end{eqnarray*}
Now, by an integration by parts in the first term of the right-hand side of the previous formula, using the fact that $A_{[\theta]}$ is divergence free and the symmetry of the operator $\mathcal{L}^{\alpha}$ we have that the expression above is equal to zero, so the quantity $\displaystyle{\int_{\mathbb{R}^n}\theta(t-s,x)\psi(s,x)dx}$, remains constant in time. We only have to set $s=0$ and $s=\frac{t}{2}$.\\

Since the drifts $A_{[(-\Delta)^{-\frac{\eta}{2}}\theta]}$ and $A_{[(-\Delta)^{\frac{\eta}{2}}\theta]}$ used in Theorems \ref{Theoreme3} and \ref{Theoreme4} are divergence free, the proof follows the same lines for these cases.  \hfill$\blacksquare$\\

This result is borrowed from \cite{KN} and was used in \cite{DCHSM}, but we have made here a slight modification which is crucial in this article. Indeed, in those previous articles, Theorem \ref{Transfert} was used to replace the H\"older-Hardy duality bracket  $\langle \theta(t,\cdot),\psi_0\rangle_{\mathcal{C}^{\lambda}\times h^{\sigma}}$ into a simpler $L^{\infty}-L^{1}$ bracket $\langle \theta_0,\psi(t,\cdot)\rangle_{L^{\infty}\times L^{1}}$ (the initial data $\theta_{0}$ was assumed in $L^{\infty}$), and then the proof of the H\"older regularity was reduced to a suitable $L^{1}$ control for deformed molecules $\psi(t, \cdot)$ since 
$$|\langle \theta_0,\psi(t,\cdot)\rangle_{L^{\infty}\times L^{1}}|\leq \|\theta_0\|_{L^\infty}\|\psi(t,\cdot)\|_{L^1}.$$

Now, with this previous theorem, we replace the bracket  $\langle \theta(t,\cdot),\psi_0\rangle_{\mathcal{C}^{\lambda}\times h^{\sigma}}$ by $\langle \theta(t/2, \cdot),\psi(t/2,\cdot)\rangle_{L^{p}\times L^{p'}}$ where $\frac{1}{p}+\frac{1}{p'}=1$ and this bracket in time $t/2$ gives us \emph{more} information than  $\langle \theta_0,\psi(t,\cdot)\rangle_{L^{p}\times L^{p'}}$ and the modification is not only related to the $L^{p}$ setting.\\

Indeed, by the maximum principle we have at time $t/2$ the inequalities
\begin{equation}\label{Dualite}
|\langle \theta(t/2, \cdot),\psi(t/2,\cdot)\rangle_{L^{p}\times L^{p'}}|\leq \|\theta(t/2,\cdot)\|_{L^{p}}\|\psi(t/2,\cdot)\|_{L^{p'}}\leq \|\theta_{0}\|_{L^{p}}\|\psi(t/2,\cdot)\|_{L^{p'}},
\end{equation}
thus we only need to control the $L^{p'}$-behavior of $\psi(t/2,\cdot)$.\\

But, at time $t/2$ we have $\theta \in L^{\infty}(L^{p})\cap L^{p}([0,T],\dot{B}^{\frac{\alpha}{p},p}_p)$ from which we deduce the uniform controls given in expression (\ref{ControlUniformes}) for the drift: in particular at time $t/2$, due to the \textcolor{black}{boundedness} properties of the drift $A_{[\theta]}$ we have a Besov control for the velocity field which \emph{was not available} at time $t=0$ and this additional information is the key of our computations. This is a specificity of the nonlinear case.
\begin{Remarque}
Note that the Besov control hypothesis is the same for Theorems \ref{Theoreme1}-\ref{Theoreme4} and what really changes from one theorem to another is the Morrey-Campanato boundedness properties. 
Since we have in these cases the uniform controls stated in Remark \ref{RemarqueControlUniformeTransport}, the corresponding smoothness degree $\alpha_{0}$ is directly given by the relationship (\ref{RelationMorrey}). We thus only have to study how to use the Besov information to obtain a new smoothness degree $\alpha$ in order to establish the comparison between the Besov and the Morrey-Campanato regularity. 
\end{Remarque}

\textbf{Proof of Theorem \ref{Theoreme1}.} From the transfer property given in Theorem \ref{Transfert} and with the inequalities (\ref{Dualite}) we only need to study the control of the $L^{p'}$ norm of $\psi(t,\cdot)$ and due to the maximum principle we can divide our proof into two steps following the molecule's size. Indeed, if $r\geq 1$ we have
$$\|\theta_0\|_{L^p}\|\psi(t/2,\cdot)\|_{L^{p'}}\leq \|\theta_0\|_{L^p}\|\psi_0\|_{L^{p'}}\leq C r^{-n+\frac{n}{p'}-\gamma} \|\theta_0\|_{L^p},$$
but, since $r\geq 1$, we obtain $|\langle \theta(t,\cdot), \psi_0\rangle_{\mathcal{C}^{\gamma}\times h^{\sigma}}|<+\infty$ for all big molecules.\\

It remains to treat the $L^{p'}$ control for \textit{small} molecules and this is done in the following theorem:
\begin{Theoreme}\label{TheoL1control}
Let $\psi_{0}$ be a small molecule and consider $\psi(t/2,\cdot)$ the associated solution of the backward problem (\ref{Evolution01}). If we assume moreover $n(n-1)\vee 2n<p<+\infty$, then there exists $C>0$ such that for a given time $0<T_0<t/2$ we have the following control of the $L^{p'}$-norm.
$$\|\psi(t/2,\cdot)\|_{L^{p'}}\leq C T_0^{-n+\frac{n}{p'}-\gamma}\qquad (T_0<t/2<T),$$
where $0<\gamma<\frac{1}{4}$. 
\end{Theoreme}
Accepting for a while this result, we have then a good control over the quantity $\|\psi(t/2,\cdot)\|_{L^{p'}}$  for all $0<r<1$ and getting back to (\ref{Dualite}) we obtain that $|\langle \theta(t,\cdot), \psi_0\rangle_{\mathcal{C}^{\gamma}\times h^{\sigma}}|$ is always bounded for $T_0<t/2<T$ and for any molecule $\psi_0$: we have proven Theorem \ref{Theoreme1} by a duality argument. The proof of Theorems \ref{Theoreme2}-\ref{Theoreme4} is essentially the same. \hfill $\blacksquare$

\textcolor{black}{
\begin{Remarque}[H\"older modulus of continuity]
Observe from Theorem \ref{TheoL1control} and the previous arguments that we actually have a control on the H\"older modulus of continuity of the solution $\theta $. Namely, for all $0<T_0<\frac t2< T$, $0<\gamma<\frac 14 $,
$$\|\theta(t,\cdot)\|_{{\mathcal C}^\gamma}\le CT_0^{-\frac np-\gamma}\|\theta_0\|_{L^p}.$$
Let us quote the related result obtained by \cite{CCZV} \textcolor{black}{in the particular setting of the SQG equation} who derive controls of the following type:
$$\|\theta(t,\cdot)\|_{{\mathcal C}^\gamma}\le C\|\theta_0\|_{L^\infty},$$
with $\gamma:=\gamma(\|\theta_0\|_{L^\infty})<1/4$.
\end{Remarque}
}
\subsection{Molecule's evolution}\label{SecEvolMol1}
The following theorem shows how the molecular properties are deformed with the evolution for a small time $s_0$.
\begin{Theoreme}\label{SmallGeneralisacion} Fix $t\in ]0,T[$. Set $\sigma$, $\gamma$ and $\omega$ real numbers such that $\frac{n}{n+1}<\sigma<1$, $\gamma=n(\frac{1}{\sigma}-1)$, $0<\gamma<\omega<1 $. Let $\psi(s_0, x)$ be a solution at time $s_0$ of problem (\ref{Evolution01}). We assume as for Theorem \ref{TheoL1control} that $0<\alpha<2$, $n(n-1)\vee 2n<p<+\infty$ and that the drift satisfies the hypotheses of Theorems \ref{Theoreme1}-\ref{Theoreme4}. Then, there exist positive constants $K$ and $\epsilon $ small  enough such that 
if $\psi_0$ is a small  molecule in the sense of Definition \ref{DefMolecules} for the local Hardy space $h^\sigma(\mathbb{R}^n)$, then for all time $0< s_0 \leq\epsilon r^\alpha$, we have the following estimates:
\begin{eqnarray}
\int_{\mathbb{R}^n}|\psi(s_0,x)||x-x(s_0)|^\omega dx &\leq &\big((\zeta r)^\alpha+Ks_0\big)^{\frac{\omega-\gamma}{\alpha}}  \label{SmallConcentration},\\
\|\psi(s_0,\cdot)\|_{L^\infty}&\leq & \frac{1}{\big((\zeta r)^\alpha+Ks_0\big)^{\frac{n+\gamma}{\alpha}}}\label{SmallLinftyevolution},\\
\|\psi(s_0,\cdot)\|_{L^1} &\leq & \frac{2v_n^{\frac{\omega}{n+\omega}}}{\big((\zeta r)^\alpha+Ks_0\big)^{\frac{\gamma}{ \alpha}}},\label{SmallL1evolution}
\end{eqnarray}
where $v_n$ denotes the volume of the $n$-dimensional unit ball.\\ 

The new molecule's center $x(s_0)$ used in formula (\ref{SmallConcentration}) is given by the evolution of the differential system
\begin{equation}\label{Defpointx_0Nonlinear}
\begin{cases}
x'(s)&= \overline{A_{[\theta,\rho]}}(t-s,x(s))=\displaystyle{\int_{\mathbb{R}^{n}}}A_{[\theta]}(t-s, x(s)-y)\varphi_{\rho}(y)dy,  s\in [0,s_0],\\[2mm]
x(0)&= x_0,
\end{cases}
\end{equation}
where $\varphi_{\rho}(s, x)=\frac{1}{\rho^{n}}\varphi\big(\textcolor{black}{\frac x \rho} 
\big)$, here $\varphi$ is a positive function in $\mathcal{C}^{\infty}_{0}(\mathbb{R}^{n})$ such that $supp(\varphi)\subset B(0,1)$ and where $\rho=\zeta^{\beta_1} r$ and $\beta_1>1$ will be specified later on.   
In the previous controls and in the dynamics for the evolution of the center, the  parameter $\zeta=\zeta(\alpha,\omega,\gamma,\mu)>1$, to be chosen further, is the same as in Definition \ref{DefMolecules}.
\end{Theoreme}
\begin{Remarque}\label{Remark4}
\begin{itemize}
\item[]
\item[1)] The expression (\ref{Defpointx_0Nonlinear}) explains the transport of the center of the molecule using the velocity $A_{[\theta]}$. The same definition is applicable to the drifts $A_{[(-\Delta)^{-\frac{\eta}{2}}\theta]}$ and $A_{[(-\Delta)^{\frac{\eta}{2}}\theta]}$ given in Theorems \ref{Theoreme3} and \ref{Theoreme4} respectively. We denote the new center of the molecule by $x(s_0)=x_{s_0} $.
\item[2)] By the maximum principle, it is enough to treat the case $0<((\zeta r)^\alpha+Ks_0)<1$ since $s_0$ is very small: otherwise the $L^1$ control will be trivial for time $s_0$ and beyond.
\item[3)] We will always assume $0<\gamma<\omega<\frac{1}{4}$. Moreover for Theorems \ref{Theoreme1}-\ref{Theoreme4} we assume the following relationship: $1/2<\delta<\alpha<2$.
\item[4)] The hypothesis $\alpha=\frac{p+n}{p+1}$ (which is common to Theorems \ref{Theoreme1}-\ref{Theoreme2}) can be rewritten as $\alpha=1-\frac{\alpha-n}{p}$ and this last formula will be helpful in the computations below. 
\end{itemize}
\end{Remarque}
\textbf{\textit{Proof of the Theorem \ref{SmallGeneralisacion}.}} 
We will adopt the following strategy: the first step is given by a study of the small Concentration condition (\ref{SmallConcentration}) and in a second step we prove the Height condition (\ref{SmallLinftyevolution}). With these two conditions at hand, the $L^1$ estimate (\ref{SmallL1evolution}) will be easily obtained. 
\subsubsection*{1) Concentration condition} 
For $s\in [0,s_0]$, we consider the functions $\Omega_{s}(x)=|x-x(s)|^{\omega}$ and $\psi(x)=\psi_+(x)-\psi_-(x)$ where the functions $\psi_{\pm}(x)\geq 0$ have disjoint support. We will denote by $\psi_\pm(s_0,x)$ two solutions of the dual equation (\ref{Evolution01}) at time $s_0$ with $\psi_\pm(0,x)=\psi_\pm(x)$. We can use the positivity principle (see \cite{DCHSM} for a proof, recall that for $\psi$ we are in a linear framework) and thus by linearity we have
$$|\psi(s_0,x)|=|\psi_+(s_0,x)-\psi_-(s_0,x)|\leq \psi_+(s_0,x)+\psi_-(s_0,x),$$ 
and we can write $\displaystyle{\int_{\mathbb{R}^n}|\psi(s_0,x)|\Omega_{s_0}(x)dx\leq\int_{\mathbb{R}^n}\psi_+(s_0,x)\Omega_{s_0}(x)dx+\int_{\mathbb{R}^n}\psi_-(s_0,x)\Omega_{s_0}(x)dx}$, so we only have to treat one of the integrals on the right hand side above. We have for all $s\in[0,s_0] $:
\begin{eqnarray*}
I_s&=&\left|\partial_{s} \int_{\mathbb{R}^n}\Omega_{s}(x)\psi_+(s,x)dx\right|\\
&=&\left|\int_{\mathbb{R}^n}\partial_{s} \Omega_{s}(x)\psi_+(s,x)+\Omega_{s}(x)\left(-\nabla\cdot\big[A_{[\theta]}(t-s,x)\, \psi_+(s,x)\big]-\mathcal{L}^{\alpha}\psi_+(s,x)\right)dx\right|\\
&=&\left|\int_{\mathbb{R}^n}-\nabla\Omega_{s}(x)\cdot x'(s)\psi_+(s,x)+\Omega_{s}(x)\left(-\nabla\cdot\big[A_{[\theta]}(t-s,x)\, \psi_+(s,x)\big]-\mathcal{L}^{\alpha}\psi_+(s,x)\right)dx\right|.
\end{eqnarray*}
\textcolor{black}{Integrating by parts we obtain}
\begin{equation*}
I_s=\left|\int_{\mathbb{R}^n}-\nabla\Omega_{s}(x)\cdot(x'(s)-A_{[\theta]}(t-s,x))\psi_+(s,x)-\Omega_{s}(x)\mathcal{L}^{\alpha}\psi_+(s,x)dx\right|.
\end{equation*}
Since the operator $\mathcal{L}^{\alpha}$ is symmetric and using the definition of $x'(s)$ given in (\ref{Defpointx_0Nonlinear}) we have
\begin{eqnarray}
I_s & \leq& c\int_{\mathbb{R}^n}|x-x(s)|^{\omega-1}|A_{[\theta]}(t-s,x)-\overline{A_{[\theta,\rho]}}(t-s,x(s))| |\psi_+(s,x)|dx + c\int_{\mathbb{R}^n}\big|\mathcal{L}^{\alpha}\Omega_{s}(x)\big|\, |\psi_+(s,x)|dx\nonumber\\[3mm]
I_{s} & \leq& c(I_{s,1}+I_{s,2}).\label{smallEstrella}
\end{eqnarray}
We will study separately each of the integrals $I_{s,1}$ and $I_{s,2}$ with two propositions.
\begin{Proposition}\label{PropositionBesovMorrey1} 
If $0<\gamma<\omega<\frac{1}{4}$, $0<\delta <\alpha<2$ and if $n(n-1)\vee 2n<p<+\infty$, we have for the integral $I_{s,1}$ in \eqref{smallEstrella} the estimate:
\begin{eqnarray*}
\int_0^{s_0} I_{s,1}ds&\leq &\eta_{1}\; r^{\omega-\alpha-\gamma},
\end{eqnarray*}
where $\eta_{1}=\eta_{1}(\zeta, \alpha, q,n)$ is such that the quantity $\frac{\eta_{1}}{\zeta^{(\omega-\alpha-\gamma)}}$ can be made very small and where $\zeta>1$ is the same parameter of Definition \ref{DefMolecules}. 
\end{Proposition}
\textbf{\textit{Proof.}}  Consider $0<\beta_0<1<\beta_1<2$ two parameters. Let $\zeta\gg1$ and $0<r\ll1$ be such that $\rho=\zeta^{\beta_1}r<1$. For a given $s\in [0,s_0] $, we write $\mathbb{R}^n=B_{\rho,s}\cup \bigcup_{k\geq 1}E_{k,s}$ where
\begin{equation}\label{SmallDecoupage}
B_{\rho,s}=\{x\in \mathbb{R}^n: |x-x(s)|\leq \rho\} \quad \mbox{and}\quad E_{k,s}= \{x\in \mathbb{R}^n: 2^{k-1} \rho<|x-x(s)|\leq  2^{k}\rho\},
\end{equation}
and we write
\begin{eqnarray}
I_{s,1}&=&\int_{B_{\rho,s}}|x-x(s)|^{\omega-1}|A_{[\theta]}(t-s,x)-\overline{A_{[\theta,\rho]}}(t-s,x(s))| |\psi(s,x)|dx\label{Boule1}\\
&&+\sum_{k\geq 1} \int_{E_{k,s}}|x-x(s)|^{\omega-1}|A_{[\theta]}(t-s,x)-\overline{A_{[\theta,\rho]}}(t-s,x(s))| |\psi(s,x)|dx.\label{Couronne1}
\end{eqnarray}
We will study each of these terms separately. 
\begin{enumerate}
\item[$\bullet$] For (\ref{Boule1}) we define $\rho_{0}=\zeta^{\beta_0}r$. Since $\beta_{0}<\beta_{1}$ and $\zeta>1$, we have $\rho_{0}<\rho$ and  we can then consider $B_{\rho,s}=B_{\rho_{0},s}\cup \mathcal{C}(\rho_{0}, \rho,s)$
with $ \mathcal{C}(\rho_{0}, \rho,s)= \{x\in \mathbb{R}^n: \rho_{0}<|x-x(s)|\leq  \rho\}$. So, we need to study $I_{B_{\rho_{0},s}}+I_{\mathcal{C}(\rho_{0}, \rho,s)}$ where
$$I_{B_{\rho_{0}},s}=\int_{B_{\rho_{0},s}}|x-x(s)|^{\omega-1}|A_{[\theta]}(t-s,x)-\overline{A_{[\theta,\rho]}}(t-s,x(s))| |\psi(s,x)|dx,$$ 
and $$I_{\mathcal{C}(\rho_{0}, \rho,s)}=\int_{\mathcal{C}(\rho_{0}, \rho,s)}|x-x(s)|^{\omega-1}|A_{[\theta]}(t-s,x)-\overline{A_{[\theta,\rho]}}(t-s,x(s))| |\psi(s,x)|dx.$$
We set 
\begin{equation}
\label{COND_HOLDER}
\frac{1}{p}+\frac{1}{p'}=1,\ \text{with}\ 1<p'<\frac{n}{1-\omega}\ \text{and}  \ \frac{1}{m}+\frac{1}{p}+\frac{1}{z}=1\ \text{with}\ m>\frac{n}{1-\omega}.
\end{equation}
 Then, by the H\"older inequality we have:
\begin{equation*}
\begin{split}
I_{B_{\rho_{0},s}}+I_{\mathcal{C}(\rho_{0}, \rho,s)}\leq& \|\psi(s,\cdot)\|_{L^\infty}\times \| |x-x(s)|^{\omega-1} \|_{L^{p'}(B_{\rho_{0},s})}\times\|A_{[\theta]}(t-s,\cdot)-\overline{A_{[\theta,\rho]}}(t-s,\cdot)\|_{L^{p}(B_{\rho_{0},s})}\\
&+\left(\int_{\mathcal{C}(\rho_{0}, \rho,s)} |A_{[\theta]}(t-s,\cdot)-\overline{A_{[\theta,\rho]}}(t-s,x(s))|^p dx \right)^{1/p}\\
&\times\left(\int_{\mathcal{C}(\rho_{0}, \rho,s)}|x-x(s)|^{(\omega-1)m}dx\right)^{1/m}\times\|\psi(s,\cdot)\|_{L^z(\mathcal{C}(\rho_{0}, \rho,s))}.
\end{split}
\end{equation*}
Since $1<p'<\frac{n}{1-\omega}$ and $\frac{n}{1-\omega}<m$, we can write
\begin{equation*}
\begin{split}
I_{B_{\rho_{0},s}}+I_{\mathcal{C}(\rho_{0}, \rho,s)}\leq & C \|\psi(s,\cdot)\|_{L^{\infty}} \times \rho_{0}^{\omega-1+\frac{n}{p'}}\times\|A_{[\theta]}(t-s,\cdot)-\overline{A_{[\theta,\rho]}}(t-s,x(s))\|_{L^{p}(B_{\rho_{0},s})}\\
&+\left(\int_{\mathcal{C}(\rho_{0}, \rho,s)} |A_{[\theta]}(t-s,\cdot)-\overline{A_{[\theta,\rho]}}(t-s,x(s))|^p dx \right)^{1/p}\\
&\times\bigg({\rho_{0}}^{m(\omega-1)+n}-\rho^{m(\omega-1)+n}\bigg)^{1/m}\times\|\psi(s,\cdot)\|_{L^{z}(\mathcal{C}(\rho_{0}, \rho,s))}, 
\end{split}
\end{equation*}
and since $\rho_{0}<\rho$ we obtain
\begin{eqnarray}
I_{B_{\rho_{0},s}}+I_{\mathcal{C}(\rho_{0}, \rho,s)}&\leq &C  \left( \|\psi(s,\cdot)\|_{L^{\infty}} \times \rho_{0}^{\omega-1+\frac{n}{p'}}+\big({\rho_{0}}^{m(\omega-1)+n}-\rho^{m(\omega-1)+n}\big)^{1/m}\times\|\psi(s,\cdot)\|_{L^{z}(B_{\rho,s})}\right)\nonumber\\[2mm]
& &\times \|A_{[\theta]}(t-s,\cdot)-\overline{A_{[\theta,\rho]}}(t-s,x(s))\|_{L^{p}(B_{\rho,s})}.\label{EstimationAvantBesov}
\end{eqnarray}

\begin{Remarque}[Morrey-Campanato versus Besov] \label{dicho_MC_B}
It is worth noting here that it is precisely at this point that we will use the Besov norm control instead of the Morrey-Campanato boundedness information. Indeed, having obtained the inequality (\ref{EstimationAvantBesov}) we have the choice to use either the first inequality of (\ref{ControlUniformes}) to control the quantity $\|A_{[\theta]}(t-s,\cdot)-\overline{A_{[\theta,\rho]}}(t-s,x(s))\|_{L^{p}(B_{\rho,s})}$, which will lead us (following \cite{DCHSM}) to relationship (\ref{RelationMorrey}), or the second inequality of (\ref{ControlUniformes}) which will allow us to consider a different smoothness degree for the L\'evy-type operator: following the choice made, we will observe a very interesting competition between the Morrey-Campanato regularity versus the Besov regularity.\\
\end{Remarque}
We have now the following lemma, proven in the appendix, which is the key of our computations
\begin{Lemme}\label{LemBesov1} 
Under the hypotheses of Theorems \ref{Theoreme1}-\ref{Theoreme2} and if $\overline{A_{[\theta,\rho]}}(t-s,x(s))$ is given by (\ref{Defpointx_0Nonlinear}), we have for $1<p<+\infty$,  $0<\alpha<2$ and almost all $s\in [0,T] $:
$$\|A_{[\theta]}(t-s,\cdot)-\overline{A_{[\theta,\rho]}}(t-s,x(s))\|_{L^p(B_{\rho,s})}\leq C \rho^{\frac{\alpha}{p}}\|A_{[\theta]}(t-s,\cdot)\|_{\dot{B}^{\alpha/p,p}_{p}(\mathbb{R}^{n})}.\\[5mm]$$ 
\end{Lemme}

\begin{Remarque}
This key lemma admits the following generalization that will be used in Theorems \ref{Theoreme3} and \ref{Theoreme4}. Indeed, if $0<\eta<1$ is a fixed (small) parameter such that respectively $0<\frac{\alpha}{p}+\eta<1$ and $0<\frac{\alpha}{p}-\eta<1$, then for almost all $s\in [0,T]$:
\begin{eqnarray}
\|A_{[(-\Delta)^{-\frac{\eta}{2}}\theta]}(t-s,\cdot)-\overline{A_{[(-\Delta)^{-\frac{\eta}{2}}\theta,\rho]}}(t-s,x(s))\|_{L^p(B_{\rho,s})}& \leq &C \rho^{\frac{\alpha}{p}+\eta}\label{EquationKeyMoins}\\ 
&\times &\|A_{[(-\Delta)^{-\frac{\eta}{2}}\theta]}(t-s,\cdot)\|_{\dot{B}^{\alpha/p+\eta,p}_{p}(\mathbb{R}^{n})}\nonumber
\end{eqnarray}
\begin{eqnarray}
\|A_{[(-\Delta)^{\frac{\eta}{2}}\theta]}(t-s,\cdot)-\overline{A_{[(\Delta)^{\frac{\eta}{2}}\theta,\rho]}}(t-s,\cdot)\|_{L^p(B_{\rho,s})}&\leq &C \rho^{\frac{\alpha}{p}-\eta}\label{EquationKeyPlus}\\
&\times &\|A_{[(-\Delta)^{\frac{\eta}{2}}\theta]}(t-s,x(s))\|_{\dot{B}^{\alpha/p-\eta,p}_{p}(\mathbb{R}^{n})}.\nonumber
\end{eqnarray}
The proof of these inequalities follows the same lines as of the one presented in Lemma \ref{LemBesov1}.
\end{Remarque}

Now we apply Lemma \ref{LemBesov1} in (\ref{EstimationAvantBesov}) and we obtain
\begin{eqnarray*}
I_{B_{\rho_{0},s}}+I_{\mathcal{C}(\rho_{0}, \rho,s)}&\leq &C  \left( \|\psi(s,\cdot)\|_{L^{\infty}} \times \rho_{0}^{\omega-1+\frac{n}{p'}}+\big({\rho_{0}}^{m(\omega-1)+n}-\rho^{m(\omega-1)+n}\big)^{1/m}\times\|\psi(s,\cdot)\|_{L^{z}(B_{\rho,s})}\right)\nonumber\\[2mm]
& &\times \rho^{\frac{\alpha}{p}} \|A_{[\theta]}(t-s,\cdot)\|_{\dot{B}^{\alpha/p,p}_{p}(\mathbb{R}^{n})} .
\end{eqnarray*} 
If we define
\begin{equation}\label{DefEpsilon}
\varepsilon=\frac{\ln\big[1-\zeta^{(\beta_1-\beta_0)(m(\omega-1)+n)}\big]}{(m (\omega-1)+n) \beta_0\ln(\zeta)},
\end{equation}
which is a positive quantity since $m(\omega-1)+n<0$ and recalling that $\rho=\zeta^{\beta_1}r $ and $\rho_{0}=\zeta^{\beta_0}r $, \textcolor{black}{with $\beta_1>\beta_0 $ and $\zeta>1 $}, we obtain
$$\left( \rho_{0}^{m(\omega-1)+n}-\rho^{m(\omega-1)+n}\right)^{1/m}=(\zeta^{\beta_0(1+\varepsilon)} r)^{\omega-1+\frac{n}{m}},$$
and we can write
\begin{eqnarray*}
I_{B_{\rho_{0},s}}+I_{\mathcal{C}(\rho_{0}, \rho,s)}& \leq  & C\bigg((\zeta^{\beta_0}r)^{\omega-1+\frac{n}{p'}}\times \|\psi(s,\cdot)\|_{L^{\infty}} +(\zeta^{\beta_0(1+\varepsilon)} r)^{\omega-1+\frac{n}{m}}\times \|\psi(s,\cdot)\|_{L^{z}(B_\rho,s)}\bigg)\\
& &\times (\zeta^{\beta_1}r)^{\frac{\alpha}{p}} \times \|A_{[\theta]}(t-s,\cdot)\|_{\dot{B}^{\alpha/p,p}_{p}(\mathbb{R}^{n})}.  
\end{eqnarray*} 
At this point, we use the maximum principle for the molecule $\psi(s,\cdot)$ \textcolor{black}{(as well as \eqref{Hipo2} and Remark \ref{Remark2}) to derive}:
$$\|\psi(s,\cdot)\|_{L^{\infty}}\leq \|\psi_{0}\|_{L^{\infty}}\leq (\zeta r)^{-(n+\gamma)}\hspace{20mm} \|\psi(s,\cdot)\|_{L^{z}(B_{\rho,s})}\leq \|\psi_{0}\|_{L^{z}(\mathbb{R}^{n})}\leq C(\zeta r)^{-n+\frac{n}{z}-\gamma},$$
which lead us to 
\begin{eqnarray*}
I_{B_{\rho_{0},s}}+I_{\mathcal{C}(\rho_{0}, \rho,s)}& \leq  & C\bigg((\zeta^{\beta_0}r)^{\omega-1+\frac{n}{p'}}\times (\zeta r)^{-(n+\gamma)} +(\zeta^{\beta_0(1+\varepsilon)} r)^{\omega-1+\frac{n}{m}}\times (\zeta r)^{-n+\frac{n}{z}-\gamma}\bigg)\\
& &\times (\zeta^{\beta_1}r)^{\frac{\alpha}{p}} \times \|A_{[\theta]}(t-s,\cdot)\|_{\dot{B}^{\alpha/p,p}_{p}(\mathbb{R}^{n})}. 
\end{eqnarray*} 
Since we have the Besov stability property for the drift given in (\ref{ControlUniformes}) -which is a common hypothesis for Theorems \ref{Theoreme1}-\ref{Theoreme2}- we can write
\begin{eqnarray*}
\int_0^{s_0} \|A_{[\theta]}(t-s,\cdot)\|_{\dot{B}^{\alpha/p,p}_{p}({\mathbb R}^n)}ds\leq  s_0^{\textcolor{black}{1/p'}}\Big(\int_0^{s_0} \|A_{[\theta]}(t-s,\cdot)\|_{\dot{B}^{\alpha/p,p}_{p}({\mathbb R}^n)}^pds\Big)^{1/p}\\
\le  \|A_{[\theta]}\|_{L^p([0,T],\dot{B}^{\alpha/p,p}_{p})} \le D_{A}\|\theta\|_{L^{p}([0,T],\dot{B}^{\alpha/p,p}_{p})}\underset{\textcolor{black}{\eqref{CTR_BESOV_LP}}}{\leq} C_{A}\|\theta_{0}\|_{L^{p}}=\mu<+\infty, 
\end{eqnarray*} 
and we obtain
\begin{eqnarray}
\int_0^{s_0}(I_{B_{\rho_{0},s}}+I_{\mathcal{C}(\rho_{0}, \rho,s)})ds & \leq & C\mu \bigg((\zeta^{\beta_0}r)^{\omega-1+\frac{n}{p'}} (\zeta r)^{-(n+\gamma)} +(\zeta^{\beta_0(1+\varepsilon)} r)^{\omega-1+\frac{n}{m}} (\zeta r)^{-n+\frac{n}{z}-\gamma}\bigg)\notag\\
&\times &(\zeta^{\beta_1}r)^{\frac{\alpha}{p}}. \label{EstimationDriftBesov1}
\end{eqnarray} 
\begin{Remarque}\label{Remarque5} In the case of Theorem \ref{Theoreme3}, we use (\ref{EquationKeyMoins}) to obtain:
\begin{eqnarray*}
\int_0^{s_0}(I_{B_{\rho_{0},s}}+I_{\mathcal{C}(\rho_{0}, \rho,s)})ds& \leq  & C\bigg((\zeta^{\beta_0}r)^{\omega-1+\frac{n}{p'}}\times (\zeta r)^{-(n+\gamma)} +(\zeta^{\beta_0(1+\varepsilon)} r)^{\omega-1+\frac{n}{m}}\times (\zeta r)^{-n+\frac{n}{z}-\gamma}\bigg)\\
& &\times (\zeta^{\beta_1}r)^{\frac{\alpha}{p}+\eta} \times \textcolor{black}{\Big(}\int_0^{s_0} \|A_{[(-\Delta)^{-\frac{\eta}{2}}\theta]}(t-s,\cdot)\|_{\dot{B}^{\alpha/p+\eta,p}_{p}(\R^n)}^{\textcolor{black}{p}}ds\textcolor{black}{\Big)^{1/p}},
\end{eqnarray*} 
moreover, using hypothesis (\ref{HypoMorreyBorne3}), Besov spaces properties, the uniform control over the Besov norm stated in 
\textcolor{black}{assumption $(b_3) $ of Theorem \ref{Theoreme3}} and the maximum principle, we have: 
\begin{eqnarray*}
\|A_{[(-\Delta)^{-\frac{\eta}{2}}\theta]}\|_{L^{p}([0,T],\dot{B}^{\alpha/p+\eta,p}_{p})}&\leq &\|(-\Delta)^{-\frac{\eta}{2}}\theta\|_{L^{p}([0,T],\dot{B}^{\alpha/p+\eta,p}_{p})}\simeq\|\theta\|_{L^{p}([0,T],\dot{B}^{\alpha/p,p}_{p})}\\
&\leq & D_{A}\| \theta_{0}\|_{L^{p}}=\mu<+\infty,
\end{eqnarray*}
and thus inequality (\ref{EstimationDriftBesov1}) may be replaced by 
\begin{eqnarray*}
\int_0^{s_0}(I_{B_{\rho_{0},s}}+I_{\mathcal{C}(\rho_{0}, \rho,s)})ds& \leq & C\mu \bigg((\zeta^{\beta_0}r)^{\omega-1+\frac{n}{p'}} (\zeta r)^{-(n+\gamma)} +(\zeta^{\beta_0(1+\varepsilon)} r)^{\omega-1+\frac{n}{m}} (\zeta r)^{-n+\frac{n}{z}-\gamma}\bigg)\\
& &\times (\zeta^{\beta_1}r)^{\frac{\alpha}{p}+\eta} 
\end{eqnarray*} 

For Theorem \ref{Theoreme4}, in a completely similar fashion, using (\ref{EquationKeyPlus}) we obtain
\begin{eqnarray*}
\int_0^{s_0}(I_{B_{\rho_{0},s}}+I_{\mathcal{C}(\rho_{0}, \rho,s)})ds &\leq & C\mu \bigg((\zeta^{\beta_0}r)^{\omega-1+\frac{n}{p'}} (\zeta r)^{-(n+\gamma)} +(\zeta^{\beta_0(1+\varepsilon)} r)^{\omega-1+\frac{n}{m}} (\zeta r)^{-n+\frac{n}{z}-\gamma}\bigg)\\
& &\times (\zeta^{\beta_1}r)^{\frac{\alpha}{p}-\eta}.\\[4mm]
\end{eqnarray*} 
\end{Remarque}

Once we have obtained inequality (\ref{EstimationDriftBesov1}) we continue our study of the molecule's evolution. Indeed, using the identity  $\alpha=1-\frac{\alpha-n}{p}$, which is common to Theorems \ref{Theoreme1}-\ref{Theoreme2} (see point 4) in Remark \ref{Remark4}) we obtain recalling as well \eqref{COND_HOLDER}:
\begin{eqnarray}
\int_0^{s_0 }(I_{B_{\rho_{0},s}}+I_{\mathcal{C}(\rho_{0}, \rho,s)})ds&\leq & C \mu\times r^{\omega-\alpha-\gamma}\nonumber\\
&  \times &\left(\zeta^{\beta_{0}(\omega-1+\frac{n}{p'})+\beta_{1}\frac{\alpha}{p}-n-\gamma}+\zeta^{\beta_0(1+\varepsilon)(\omega-1+\frac{n}{m})+\beta_{1}\frac{\alpha}{p}-n+\frac{n}{z}-\gamma}\right).\label{EstimationCoefficientsConcentration1}
\end{eqnarray} 
\begin{Remarque}\label{RemarqueRelationAlpha}
For Theorem \ref{Theoreme3}, using relationship (\ref{DegreLevy11}) we have $\alpha=1+\frac{n-\alpha}{p}-\eta$. Inequality (\ref{EstimationCoefficientsConcentration1}) becomes
\begin{eqnarray*}
\int_0^{s_0}(I_{B_{\rho_{0,s}}}+I_{\mathcal{C}(\rho_{0}, \rho,s)})ds&\leq & C \mu\times r^{\omega-\alpha-\gamma}\times \\
& &\left(\zeta^{\beta_{0}(\omega-1+\frac{n}{p'})+\beta_{1}(\frac{\alpha}{p}+\eta)-n-\gamma}+\zeta^{\beta_0(1+\varepsilon)(\omega-1+\frac{n}{m})+\beta_{1}(\frac{\alpha}{p}+\eta)-n+\frac{n}{z}-\gamma}\right).
\end{eqnarray*}
Whereas for Theorem \ref{Theoreme4}, with relationship (\ref{DegreLevy12}) we have $\alpha=1+\frac{n-\alpha}{p}+\eta$ and thus inequality (\ref{EstimationCoefficientsConcentration1}) becomes
$$I_{B_{\rho_{0},s}}+I_{\mathcal{C}(\rho_{0}, \rho,s)}\leq  C \mu\times r^{\omega-\alpha-\gamma}\times \left(\zeta^{\beta_{0}(\omega-1+\frac{n}{p'})+\beta_{1}(\frac{\alpha}{p}-\eta)-n-\gamma}+\zeta^{\beta_0(1+\varepsilon)(\omega-1+\frac{n}{m})+\beta_{1}(\frac{\alpha}{p}-\eta)-n+\frac{n}{z}-\gamma}\right).$$
\end{Remarque}
\quad\\

Once we have inequality (\ref{EstimationCoefficientsConcentration1}), we want to keep the information of the molecule's relative size $r^{\omega-\alpha-\gamma}$ and we want to absorb the constants, thus if we denote $\eta_{0,1}$ by
$$\eta_{0,1}=C \mu\left(\zeta^{\beta_{0}(\omega-1+\frac{n}{p'})+\beta_{1}\frac{\alpha}{p}-n-\gamma}+\zeta^{\beta_0(1+\varepsilon)(\omega-1+\frac{n}{m})+\beta_{1}\frac{\alpha}{p}-n+\frac{n}{z}-\gamma}\right),$$
we can write
\begin{equation}\label{Boule}
\int_0^{s_0}(I_{B_{\rho_{0,s}}}+I_{\mathcal{C}(\rho_{0}, \rho,s)}) ds\leq \eta_{0,1} r^{\omega-\alpha-\gamma}.
\end{equation}
We remark now that $\frac{\eta_{0,1}}{\zeta^{(\omega-\alpha-\gamma)}}$ can be made small, indeed we have
\begin{eqnarray}
\frac{\eta_{0,1}}{\zeta^{(\omega-\alpha-\gamma)}}&=&C \mu\left(\zeta^{\beta_{0}(\omega-1+\frac{n}{p'})+\beta_{1}\frac{\alpha}{p}-n-\gamma-(\omega-\alpha-\gamma)}+\zeta^{\beta_0(1+\varepsilon)(\omega-1+\frac{n}{m})+\beta_{1}\frac{\alpha}{p}-n+\frac{n}{z}-\gamma-(\omega-\alpha-\gamma)}\right)\nonumber\\
&=& C \mu\left(\zeta^{\beta_{0}(\omega-1+\frac{n}{p'})+\beta_{1}\frac{\alpha}{p}-n-\omega+\alpha}+\zeta^{\beta_0(1+\varepsilon)(\omega-1+\frac{n}{m})+\beta_{1}\frac{\alpha}{p}-n+\frac{n}{z}-\omega+\alpha}\right),\label{Negativity}
\end{eqnarray}
and since the parameter $\zeta\gg1$, it is enough to verify that the quantities $\beta_{0}(\omega-1+\frac{n}{p'})+\beta_{1}\frac{\alpha}{p}-n-\omega+\alpha$ and $\beta_0(1+\varepsilon)(\omega-1+\frac{n}{m})+\beta_{1}\frac{\alpha}{p}-n+\frac{n}{z}-\omega+\alpha$ can be chosen to be both negative.\\ 

If we set $\beta_{0}=1-\nu_{0}$ and $\beta_{1}=1+\nu_{1}$ with $0<\nu_{1}\ll \nu_{0}<1$ and recalling the relationship $\alpha=1-\frac{\alpha-n}{p}$, the negativity of $\beta_{0}(\omega-1+\frac{n}{p'})+\beta_{1}\frac{\alpha}{p}-n-\omega+\alpha$ is equivalent to 
\begin{equation}\label{EstimationCoefficientsConcentration2}
\nu_{1}\frac{\alpha}{p}<\nu_{0}(\frac{n}{p'}+\omega-1),
\end{equation}
but since from \eqref{COND_HOLDER},  $1<p'<\frac{n}{1-\omega}$ and $\nu_{1}\ll \nu_{0}$, we obtain that the first power of $\zeta$ in (\ref{Negativity}) is negative.\\

\begin{Remarque}
For Theorem \ref{Theoreme3}, we need to study the negativity of the quantity $\beta_{0}(\omega-1+\frac{n}{p'})+\beta_{1}(\frac{\alpha}{p}+\eta)-n-\omega+\alpha$, 
and the condition (\ref{EstimationCoefficientsConcentration2}) changes as follows
$\nu_{1}\left(\frac{\alpha}{p}+\eta\right)<\nu_{0}(\frac{n}{p'}+\omega-1)$, while for Theorem \ref{Theoreme4} we need to study $\beta_{0}(\omega-1+\frac{n}{p'})+\beta_{1}(\frac{\alpha}{p}-\eta)-n-\omega+\alpha$, and the condition obtained is $\nu_{1}\left(\frac{\alpha}{p}-\eta\right)<\nu_{0}(\frac{n}{p'}+\omega-1)$. These conditions are fulfilled as long as $1<p'<\frac{n}{1-\omega}$ and $\nu_{1}\ll \nu_{0}$.
\end{Remarque}

With the same relationships as above for $\beta_{0}$, $\beta_{1}$, $\nu_{0}$ and $\nu_{1}$, the negativity of $\beta_0(1+\varepsilon)(\omega-1+\frac{n}{m})+\beta_{1}\frac{\alpha}{p}-n+\frac{n}{z}-\omega+\alpha$ is equivalent to 
\begin{equation}\label{EstimationCoefficientsConcentration3}
(\omega-1+\frac{n}{m})(\varepsilon-\nu_{0}(1+\varepsilon))+\nu_{1}(1-\alpha+\frac{n}{p})<0,
\end{equation}
where $\frac{n}{1-\omega}<m<+\infty$ and $\varepsilon$ is given by (\ref{DefEpsilon}). Since $\varepsilon=\varepsilon(\nu_{0}, \nu_{1})\underset{\nu_{0}\to 0}{\longrightarrow}+\infty$, the negativity of this exponent of $\zeta$ follows for $\nu_{0}$ small enough. 

\begin{Remarque} As before, for Theorem \ref{Theoreme3} we need to study the quantity
$\beta_0(1+\varepsilon)(\omega-1+\frac{n}{m})+\beta_{1}(\frac{\alpha}{p}+\eta)-n+\frac{n}{z}-\omega+\alpha$, but since we have $\alpha=1+\frac{n-\alpha}{p}-\eta$ (see Remark \ref{RemarqueRelationAlpha}), the negativity of this quantity is equivalent to the condition (\ref{EstimationCoefficientsConcentration3}). For Theorem \ref{Theoreme4} we study
$\beta_0(1+\varepsilon)(\omega-1+\frac{n}{m})+\beta_{1}(\frac{\alpha}{p}-\eta)-n+\frac{n}{z}-\omega+\alpha$, but since  $\alpha=1+\frac{n-\alpha}{p}+\eta$ (see also Remark \ref{RemarqueRelationAlpha}), we obtain the same condition (\ref{EstimationCoefficientsConcentration3}).
\end{Remarque}

Thus, since all the exponents in (\ref{Negativity}) are negative, it is possible to take $\zeta$ big enough in order to make $\frac{\eta_{0,1}}{\zeta^{(\omega-\alpha-\gamma)}}$ small.\\[5mm]
\item[$\bullet$] We study now (\ref{Couronne1}). For $s\in [0,s_0] $, let us denote by $I_{E_{k,s}}$ the integral 
$$I_{E_{k,s}}=\int_{E_{k,s}}|x-x(s)|^{\omega-1}|A_{[\theta]}(t-s,\cdot)-\overline{A_{[\theta,\rho]}}(t-s,x(s))| |\psi(s,x)|dx.$$
As $\omega-1<0$, over $E_{k,s}$ we have $|x-x(s)|^{\omega-1}\leq C (2^{k}\rho)^{\omega-1}$ and we write
\begin{eqnarray*}
I_{E_{k,s}}&\leq & C(2^{k}\rho)^{\omega-1}\int_{B_{2^k \rho,s}}|A_{[\theta]}(t-s,\cdot)-\overline{A_{[\theta,\rho]}}(t-s,x(s))| |\psi(s,x)|dx,
\end{eqnarray*}
where we have denoted $B_{2^k \rho,s}=B(x(s),2^k \rho)$, then we obtain
\begin{equation}\label{InegaliteCouronnesBesov}
I_{E_{k,s}} \leq C(2^{k}\rho)^{\omega-1} \|A_{[\theta]}(t-s,\cdot)-\overline{A_{[\theta,\rho]}}(t-s,x(s))\|_{L^p(B_{2^k \rho,s})} \|\psi(s,\cdot)\|_{L^{p'}(B_{2^k \rho,s})},
\end{equation}
where we used the Hölder inequality with $\frac{1}{p}+\frac{1}{p'}=1$. We will need the following lemma which is proved in the appendix. 
\begin{Lemme}\label{LemBesov2} Under the general hypotheses of Theorems \ref{Theoreme1}-\ref{Theoreme2} and if $\overline{A_{[\theta,\rho]}}(t-s,\cdot)$ is given by (\ref{Defpointx_0Nonlinear}), we have for almost all $s\in [0,s_0] $,
$$\|A_{[\theta]}(t-s,\cdot)-\overline{A_{[\theta,\rho]}}(t-s,x(s))\|_{L^p(B_{2^k \rho,s})}\leq C(2^k \rho)^{\frac{\alpha}{p}} \times 2^{\frac{kn}{p}} \times\|A_{[\theta]}(t-s,\cdot)\|_{\dot{B}^{\alpha/p,p}_{p}(\mathbb{R}^{n})}
.$$ 
\end{Lemme}

\begin{Remarque}\label{LA_REM_COUR_DYAD_PERTURB}
This lemma admits the following generalization: if $0<\eta<1$ is a fixed (small) parameter such that $0<\frac{\alpha}{p}+\eta<1$ and $0<\frac{\alpha}{p}-\eta<1$, then we have for almost all $s\in [0,s_0] $,
\begin{eqnarray}
\|A_{[(-\Delta)^{-\frac{\eta}{2}}\theta]}(t-s,\cdot)-\overline{A_{[(-\Delta)^{-\frac{\eta}{2}}\theta,\rho]}}(t-s,\cdot)\|_{L^p(B_{2^k \rho,s})}&\leq & C (2^k \rho)^{\frac{\alpha}{p}+\eta} \label{EquationKeyMoins1}\\
&\times &2^{\frac{kn}{p}} \|A_{[(-\Delta)^{-\frac{\eta}{2}}\theta]}(t-s,\cdot)\|_{\dot{B}^{\alpha/p+\eta,p}_{p}(\mathbb{R}^{n})}.\nonumber
\end{eqnarray}
The previous inequality will be used in Theorem \ref{Theoreme3}, and for Theorem \ref{Theoreme4} we will need the following estimate, for almost all $s\in [0,s_0] $:
\begin{eqnarray}
\|A_{[(-\Delta)^{\frac{\eta}{2}}\theta]}(t-s,\cdot)-\overline{A_{[(-\Delta)^{\frac{\eta}{2}}\theta,\rho]}}(t-s,\cdot)\|_{L^p(B_{2^k \rho,s})}&\leq & C (2^k\rho)^{\frac{\alpha}{p}-\eta} \label{EquationKeyPlus1}\\
& \times &2^{\frac{kn}{p}} \|A_{[(-\Delta)^{\frac{\eta}{2}}\theta]}(t-s,\cdot)\|_{\dot{B}^{\alpha/p-\eta,p}_{p}(\mathbb{R}^{n})}.\nonumber
\end{eqnarray}
The proof of these inequalities follows the same lines as of the one presented in Lemma \ref{LemBesov2}.
\end{Remarque}
If we apply Lemma \ref{LemBesov2} in inequality (\ref{InegaliteCouronnesBesov}) we obtain
\begin{eqnarray*}
I_{E_{k,s}} &\leq & C(2^{k}\rho)^{\omega-1} \left((2^k \rho)^{\frac{\alpha}{p}} \times 2^{\frac{kn}{p}} \times\|A_{[\theta]}(t-s,\cdot)\|_{\dot{B}^{\alpha/p,p}_{p}(\mathbb{R}^{n})}\|\psi(s,\cdot)\|_{L^{p'}(\mathbb{R}^{n})}\right),
\end{eqnarray*}
and using again the maximum principle for the molecule\footnote{Recall from the maximum principle that we have $\|\psi(s,\cdot)\|_{L^{p'}(\mathbb{R}^{n})}\leq \|\psi_{0}\|_{L^{p'}(\mathbb{R}^{n})}\leq C(\zeta r)^{-\frac{n}{p}-\gamma}.$} $\psi(s,\cdot)$ we have
$$I_{E_{k,s}} \leq  C(2^{k}\rho)^{\omega-1} \times (\zeta r)^{-n+\frac{n}{p'}-\gamma}\times \left((2^k \rho)^{\frac{\alpha}{p}} \times 2^{\frac{kn}{p}} \times\|A_{[\theta]}(t-s,\cdot)\|_{\dot{B}^{\alpha/p,p}_{p}(\mathbb{R}^{n})}\right).$$
Now, in the framework of Theorems \ref{Theoreme1}-\ref{Theoreme2} we have by hypothesis the Besov \textcolor{black}{boundedness} property  $\|A_{[\theta]}\|_{L^{p}([0,T],\dot{B}^{\alpha/p,p}_{p})}\leq D_{A}\|\theta\|_{L^{p}([0,T],\dot{B}^{\alpha/p,p}_{p})}\leq D_{A}\|\theta_{0}\|_{L^{p}}=\mu<+\infty$,  and we can write recalling from \eqref{COND_HOLDER} that $\frac 1p+\frac{1}{p'}=1 $:
\begin{equation}\label{EstimationCouronnesTheo}
\int_0^{s_0} I_{E_{k,s}}ds \leq  C\mu(2^{k}\rho)^{\omega-1} \times (\zeta r)^{-\frac{n}{p}-\gamma}\times \left((2^k \rho)^{\frac{\alpha}{p}} \times 2^{\frac{kn}{p}}\right).
\end{equation}
\begin{Remarque}\label{Remarque6} For Theorem \ref{Theoreme3}, using the inequality (\ref{EquationKeyMoins1}) and the estimates given in  \textcolor{black}{Remark \ref{LA_REM_COUR_DYAD_PERTURB}} we have instead of (\ref{EstimationCouronnesTheo}) the following inequality
$$\int_0^{s_0} I_{E_{k,s}} ds\leq  C\mu(2^{k}\rho)^{\omega-1} \times (\zeta r)^{-\frac{n}{p}-\gamma}\times \left((2^k \rho)^{\frac{\alpha}{p}+\eta} \times 2^{\frac{kn}{p}}\right).$$
For Theorem \ref{Theoreme4}, using (\ref{EquationKeyPlus1}) and Remark \ref{Remarque5} we obtain
$$\int_0^{s_0}I_{E_{k,s}} ds\leq  C\mu(2^{k}\rho)^{\omega-1} \times (\zeta r)^{-\frac{n}{p}-\gamma}\times 
\left((2^k \rho)^{\frac{\alpha}{p}-\eta} \times 2^{\frac{kn}{p}}\right).$$
\end{Remarque}
\quad\\

Once we have at our disposal inequality (\ref{EstimationCouronnesTheo}), recalling that  $\rho=\zeta^{\beta_{1}}r$ and $\alpha=1-\frac{\alpha-n}{p}$, we have the following expression
\begin{eqnarray}\label{EstimationCouronnesSomme}
\int_0^{s_0} I_{E_{k,s}}ds&\leq & C\mu 2^{k(\omega-1+\frac{\alpha+n}{p})} \times\zeta^{\beta_{1}(\omega-1+\frac{\alpha}{p})-\frac{n}{p}-\gamma}\times r^{\omega-\alpha-\gamma}.
\end{eqnarray}
Since $0<\omega<\frac{1}{4}$, $0<\alpha<2$, $2\leq n$ and $n(n-1)\vee 2n<p<+\infty$, we have that  $\omega-1+\frac{\alpha+n}{p}<0$ and we can sum over $k$ in order to obtain
\begin{eqnarray}\label{Coronak}
\sum_{k\geq 1}\int_0^{s_0}I_{E_{k,s}}ds&\leq& C\mu \zeta^{\beta_{1}(\omega-1+\frac{\alpha}{p})-\frac{n}{p}-\gamma}\times r^{\omega-\alpha-\gamma}.
\end{eqnarray}
\begin{Remarque} Using Remark \ref{Remarque6}, for Theorem \ref{Theoreme3} we have instead of (\ref{EstimationCouronnesSomme}) the following inequality
$$\int_0^{s_0}I_{E_{k,s}}ds\leq  C\mu 2^{k(\omega-1+\frac{\alpha+n}{p}+\eta)} \times\zeta^{\beta_{1}(\omega-1+\frac{\alpha}{p}+\eta)-\frac{n}{p}-\gamma}\times r^{\omega-\alpha-\gamma},$$
and thus, if $\omega-1+\frac{\alpha+n}{p}+\eta<0$ (which is the case since $0<\omega<\frac{1}{4}$, $0<\alpha<2$, $2\leq n$, $n(n-1)\vee 2n<p<+\infty$ and $0<\eta<1$ is assumed small enough), summing over $k$ we obtain:
\begin{equation}\label{EquationKeyMoins2}
\sum_{k\geq 1}\int_0^{s_0}I_{E_{k,s}}ds\leq C\mu \zeta^{\beta_{1}(\omega-1+\frac{\alpha}{p}+\eta)-\frac{n}{p}-\gamma}\times r^{\omega-\alpha-\gamma}.
\end{equation}
For Theorem \ref{Theoreme4} we have $\displaystyle{\int_0^{s_0}}I_{E_{k,s}}ds\leq  C\mu 2^{k(\omega-1+\frac{\alpha+n}{p}-\eta)} \times\zeta^{\beta_{1}(\omega-1+\frac{\alpha}{p}-\eta)-\frac{n}{p}-\gamma}\times r^{\omega-\alpha-\gamma}$, but since $\omega-1+\frac{\alpha+n}{p}-\eta<0$ we obtain
\begin{equation}\label{EquationKeyPlus2}
\sum_{k\geq 1}\int_0^{s_0}I_{E_{k,s}}ds\leq C\mu \zeta^{\beta_{1}(\omega-1+\frac{\alpha}{p}-\eta)-\frac{n}{p}-\gamma}\times r^{\omega-\alpha-\gamma}.
\end{equation}
\end{Remarque}
\quad\\

Now, coming back to the expression (\ref{Coronak}), setting $\eta_{1,1}:=C\mu \zeta^{\beta_{1}(\omega-1+\frac{\alpha}{p})-\frac{n}{p}-\gamma}$, we observe that $\frac{\eta_{1,1}}{\zeta^{(\omega-\alpha-\gamma)}}$ can be made small. Indeed:
$$\frac{\eta_{1,1}}{\zeta^{(\omega-\alpha-\gamma)}}=C\mu \zeta^{\beta_{1}(\omega-1+\frac{\alpha}{p})-\frac{n}{p}-\omega+\alpha}.$$
Since $\zeta\gg1$ we only need to verify that $\beta_{1}(\omega-1+\frac{\alpha}{p})-\frac{n}{p}-\omega+\alpha<0$, but since $\beta_{1}=1+\nu_{1}$ and $\alpha=1-\frac{\alpha-n}{p}$, the negativity of this exponent is equivalent to $\omega+\frac{n}{p}<\alpha$ which is the case by hypothesis since $0<\omega<\frac{1}{4}$, $1<\alpha<2$ (for Theorems \ref{Theoreme1} \& \ref{Theoreme2}) and $n(n-1)\vee 2n<p<+\infty$ with $2\leq n$. 

\begin{Remarque} In order to absorb the constants in the expressions (\ref{EquationKeyMoins2}) and (\ref{EquationKeyPlus2}) we need to study the negativity of the quantities 
$$\beta_{1}(\omega-1+\frac{\alpha}{p}+\eta)-\frac{n}{p}-\omega+\alpha<0\quad \mbox{ and }\quad \beta_{1}(\omega-1+\frac{\alpha}{p}-\eta)-\frac{n}{p}-\omega+\alpha<0,$$
but using the fact that $\frac{\alpha}{p}=1+\frac{n}{p}-\alpha-\eta$ for Theorem \ref{Theoreme3} and $\frac{\alpha}{p}=1+\frac{n}{p}-\alpha+\eta$ for Theorem \ref{Theoreme4} (see Remark \ref{RemarqueRelationAlpha}), the negativity of these terms follows from the same condition $\omega+\frac{n}{p}<\alpha$.\\

Remark also that, in the case of Theorem \ref{Theoreme3}, the condition $\omega+\frac{n}{p}<\alpha$ gives a lower bound for the smoothness degree $\alpha$ given in (\ref{DegreLevy11}). Moreover, since $0<\omega<\frac{1}{4}$ and $n(n-1)\vee 2n<p<+\infty$, we can obtain that $1/2<\alpha<1$.\\

More accurate bounds can be obtained, but this is enough to our purposes, in particular in the case of Theorem \ref{Theoreme3} we can consider smoothness degrees in the interval $1/2<\alpha<1$. 
\end{Remarque}
\end{enumerate}

Now, in order to finish the proof of the Proposition \ref{PropositionBesovMorrey1}, it remains to gather the estimate (\ref{Boule}) with (\ref{Coronak}) to obtain
\begin{eqnarray*}
\int_0^{s_0}I_{s,1}ds&=&\int_0^{s_0}\big(I_{B_{\rho_{0},s}}+I_{\mathcal{C}(\rho_{0}, \rho,s)}+\sum_{k\geq 1}I_{E_{k},s}\big) ds\\
&\leq &r^{\omega-\alpha-\gamma}(\eta_{0,1} +\eta_{1,1})=r^{\omega-\alpha-\gamma}\eta_{1},
\end{eqnarray*}
where $\eta_{1}=\eta_{0,1}+\eta_{1,1}$ is such that $\frac{\eta_{1}}{\zeta^{(\omega-\alpha-\gamma)}}$ can be made very small. This concludes the proof of Proposition \ref{PropositionBesovMorrey1}. \hfill$\blacksquare$\\

With the several remarks stated all along the proof of Proposition \ref{PropositionBesovMorrey1}, we also have this estimate in the framework of Theorem \ref{Theoreme3} and Theorem \ref{Theoreme4}.\\

The second step is given by the following proposition that relies only in the properties of the L\'evy-type operator. 
\begin{Proposition}\label{PropositionBesovMorrey2} 
For the integral $I_{s,2}$ in (\ref{smallEstrella}) we have the following inequality:
\begin{equation*}
\int_{0}^{s_{0}}I_{s,2}ds\leq \eta_{2}r^{\omega-\alpha-\gamma},
\end{equation*}
where $\omega-\delta+\textcolor{black}{\frac n p}<0$ with $1/2<\delta <\alpha<2$ and where $\frac{\eta_{2}}{\zeta^{(\omega-\alpha-\gamma)}}$ can be made very small.
\end{Proposition}
Remark in particular that since $0<\omega<\frac{1}{4}$ and $n(n-1)\vee 2n<p<+\infty$, then the condition $\omega-\delta+\frac n{p}<0$ is fulfilled. As Proposition \ref{PropositionBesovMorrey2} only involves the properties of the L\'evy-type operator $\mathcal{L}^{\alpha}$, its proof follows the same lines given in \cite{DCHSM}, Lemma 4.8. We refer to this article for the details.\\

With the help of Propositions \ref{PropositionBesovMorrey1} and \ref{PropositionBesovMorrey2}, we go back to (\ref{smallEstrella}) and we have
\begin{eqnarray*}
\int_0^{s_0}I_{s} ds & \leq& c\int_0^{s_0}(I_{s,1}+I_{s,2})ds\leq c(\eta_{1}+\eta_{2})r^{\omega-\alpha-\gamma}.
\end{eqnarray*}
Using the initial concentration condition \eqref{Hipo2} one gets:
\begin{eqnarray*}
\int_{{\mathbb R}^n}|x-x(s_0)|^\omega |\psi(s_0,x)| dx&\leq & (\zeta r)^{\omega-\gamma}+2c(\eta_1+\eta_2)r^{\omega-\alpha-\gamma} s_0\\
&\le& \left((\zeta r)^\alpha+2c \frac{\alpha}{\omega-\gamma}(\eta_{1}+\eta_{2})  \frac{s_0}{\zeta^{\omega-\alpha-\gamma}}\right)^{\frac{\omega-\gamma}{\alpha}} =\big((\zeta r)^\alpha+K s_0\big)^{\frac{\omega-\gamma}{\alpha}},
\end{eqnarray*}
and we want to make the quantity $K=2c\frac{\alpha}{\omega-\gamma}(\eta_{1}+\eta_{2})\frac{1}{\zeta^{(\omega-\alpha-\gamma)}}$ very small. By the computations above, it is possible to choose $\zeta$ big enough in order to obtain the inequality
\begin{equation}\label{ConditionK11}
K=2c\frac{\alpha}{\omega-\gamma}(\eta_{1}+\eta_{2})\frac{1}{\zeta^{(\omega-\gamma-\alpha)}}\leq \left(\frac{\alpha}{n+\gamma} \right) \overline{c}_1 \times \mathfrak{c},
\end{equation}
where the constant $ \overline{c}_1$ is related to the hypotheses (\ref{DefKernelLevy}) assumed for the Lévy operator $\mathcal{L}$ and where the constant $\mathfrak{c}:=\frac{v_n (5^n-1)-\sqrt{2v_n}5^{n-\omega}}{2 \times 5^{n+\alpha}}$ is chosen as in the proof of Theorem 10 of \cite{DCHSM} to handle the height condition (see also the remarks below). We finally obtain the following inequality for the concentration condition with an appropriate control of the constant $K$:
$$\int_{{\mathbb R}^n}|x-x(s_0)|^\omega \psi(s_0,x) dx\leq  \big((\zeta r)^\alpha+Ks_0 \big)^{\frac{\omega-\gamma}{\alpha}}.$$
This concludes the study of the concentration condition.
\subsubsection*{2) Height condition and $L^1 $ estimate}
\textcolor{black}{The height condition only  involves} the operator ${\mathcal L}^\alpha $, independently of the drift, since we are led to investigate the time evolution of the supremum norm in space which precisely annihilates the transport term.

This contribution can readily be analyzed following the proof of Theorem 10 in \cite{DCHSM} which precisely yields estimate \eqref{SmallLinftyevolution}. \textcolor{black}{The $L^1$ estimate} of equation \eqref{SmallL1evolution} then follows from a direct optimization over the parameter $D$ writing
\begin{eqnarray*}
\int_{\mathbb{R}^n}|\psi(s_0,x)|dx&=&\int_{\{|x-x(s_0)|< D\}}|\psi(s_0,x)|dx+\int_{\{|x-x(s_0)|\geq D\}}|\psi(s_0,x)|dx\\
&\leq & v_n D^n \|\psi(s_0,\cdot)\|_{L^\infty}+D^{-\omega}\int_{\mathbb{R}^{n}}|\psi(s_0,x)||x-x(s_0)|^\omega dx.
\end{eqnarray*}
Theorem \ref{SmallGeneralisacion} is now completely proven. \hfill$\blacksquare$

\subsection{$L^{p'}$ control and Iteration}\label{SecEvolMol2}
Once we have a good control over the quantities $\|\psi(s_0,\cdot)\|_{L^1}$ and $\|\psi(s_0,\cdot)\|_{L^\infty}$ \textcolor{black}{(from \eqref{SmallL1evolution} and \eqref{SmallLinftyevolution})}, by interpolation we obtain the following bound
$$\|\psi(s_0,\cdot)\|_{L^{p'}}\leq \|\psi(s_0,\cdot)\|_{L^1}^{\frac{1}{p'}} \|\psi(s_0,\cdot)\|_{L^\infty}^{1-\frac{1}{p'}}\leq C\left[\big((\zeta r)^\alpha+Ks_0\big)^{\frac{1}{\alpha}} \right]^{-n+\frac{n}{p'}-\gamma}.$$
We thus see with Theorem \ref{SmallGeneralisacion} that it is possible to control the $L^{p'}$ norm of the molecules $\psi$ from $0$ to a small time $s_0$, and \textcolor{black}{applying inductively} the same arguments we can extend the control from time $s_0$ to time $s_{1}$ with a small increment $s_{1}-s_{0}\sim r^{\alpha}$.  Now we can see that the smallness of the parameter $r$ and of the time increments $s_0,s_1-s_0,...,s_N-s_{N-1}$ can be compensated by the number of iterations $N$: fix a small $0<r<1$ and iterate as explained before. Since each small time increment $s_0,s_1-s_0,...,s_N-s_{N-1} $ has order $\epsilon r^\alpha$, we have $s_{N}\sim N \epsilon r^\alpha$. Thus, we will stop the iterations as soon as $ N \epsilon r^\alpha\geq T_0$. \\

The number of iterations $N=N(r)$ will depend on the smallness of the molecule's size $r$, and it is enough to consider $N(r)\sim \frac{T_0}{\epsilon r^\alpha}$ in order to obtain this lower bound for $N(r)$. Proceeding this way we will obtain $\|\psi(s_N,\cdot)\|_{L^{p'}}\leq C T_0^{-n+\frac{n}{p'}-\gamma}<+\infty$, for all molecules of size $r$. Observe that once this estimate is available, for larger times it is enough to apply the maximum principle.\\

Finally, and for all $r>0$, we obtain after a time $T_0$ a $L^{p'}$ control for small molecules and we finish the proof of the Theorem \ref{TheoL1control}. \hfill$\blacksquare$
\appendix
\mysection{Appendix}
We prove in this section Lemmas \ref{LemBesov1} and \ref{LemBesov2}.\\

\textbf{\textit{Proof of Lemma \ref{LemBesov1}.}} Using the definition of $\overline{A_{[\theta,\rho]}}$ given in (\ref{Defpointx_0Nonlinear}), and for $1<p<+\infty$ we write
$$\|A_{[\theta]}(t-s,\cdot)-\overline{A_{[\theta,\rho]}}(t-s,x(s))\|_{L^p(B_\rho)}^p=\int_{B_{\rho}}\left|A_{[\theta]}(t-s,x)-\int_{\mathbb{R}^{n}}A_{[\theta]}(t-s,x(s)-y)\varphi_{\rho}(y)dy\right|^pdx $$
\begin{eqnarray}
&\le & 2^{p-1}\Big(\int_{B_{\rho}}\left|A_{[\theta]}(t-s,x)-\int_{\mathbb{R}^{n}}A_{[\theta]}(t-s,x-y)\varphi_{\rho}(y)dy\right|^pdx\nonumber\\
& &+\int_{B_{\rho}}\left|\int_{\mathbb{R}^{n}}\big(A_{[\theta]}(t-s,x(s)-y)-A_{[\theta]}(t-s,x-y)\big)\varphi_{\rho}(y)dy\right|^pdx\Big)\nonumber\\
&\le &2^{p-1}\Big(\int_{B_{\rho}} \int_{\mathbb{R}^{n}}\left|A_{[\theta]}(t-s,x)-A_{[\theta]}(t-s,x-y)\right|^p|\varphi_{\rho}(y)|dy\;dx\nonumber\\
& &+\int_{B_{\rho}}\left|\int_{\mathbb{R}^{n}}A_{[\theta]}(t-s,y)\big(\varphi_{\rho}(x-y)-\varphi_\rho(x(s)-y)\big)dy\right|^pdx\Big),\label{PREAL_BESOV}
\end{eqnarray}
using the Jensen inequality for the last bound.
Since $\displaystyle{\int_{\R^n}}  \big(\varphi_{\rho}(x-y)-\varphi_\rho(x(s)-y)\big)dy =0$ and  $ \|\varphi_{\rho}\|_{L^{1}(\mathbb{R}^{n})}=1$,  we obtain
\begin{eqnarray*}
\left|\int_{\mathbb{R}^{n}}A_{[\theta]}(t-s,y)\big(\varphi_{\rho}(x-y)-\varphi_\rho(x(s)-y)\big)dy\right|^p=\\
\left|\int_{\mathbb{R}^{n}}\big(A_{[\theta]}(t-s,y)-A_{[\theta]}(t-s,x)\big)\big(\varphi_{\rho}(x-y)-\varphi_\rho(x(s)-y)\big)dy\right|^p\\
\le \textcolor{black}{2^{p-1}} \int_{\R^n} \Big|A_{[\theta]}(t-s,y)-A_{[\theta]}(t-s,x)\Big|^p |\varphi_\rho(x-y)-\varphi_\rho(x(s)-y)|dy.
\end{eqnarray*}
Plugging this estimate in \eqref{PREAL_BESOV} yields:
\begin{eqnarray*}
\|A_{[\theta]}(t-s,\cdot)-\overline{A_{[\theta,\rho]}}(t-s,\cdot)\|_{L^p(B_\rho)}^p\leq  2^{p-1}\Big( \int_{B_{\rho}} \int_{\mathbb{R}^{n}}\left|A_{[\theta]}(t-s,x)-A_{[\theta]}(t-s,x-y)\right|^p|\varphi_{\rho}(y)|dy\;dx\\
+\textcolor{black}{2^{p-1}}\int_{B_{\rho}} \int_{\mathbb{R}^{n}}\Big|A_{[\theta]}(t-s,y)-A_{[\theta]}(t-s,x)\Big|^p \big(\varphi_\rho(x-y)+\varphi_\rho(x(s)-y)\big) dy dx\Big).
\end{eqnarray*}
Now if we consider a small fixed parameter $s$ such that $0<sp<1$ we can write
\begin{eqnarray*}
\|A_{[\theta]}(t-s,\cdot)-\overline{A_{[\theta,\rho]}}(t-s,\cdot)\|_{L^p(B_\rho)}^p\\
\leq 2^{p-1}\Big( \int_{B_{\rho}} \int_{\mathbb{R}^{n}}\frac{|A_{[\theta]}(t-s,x)-A_{[\theta]}(t-s,x-y)|^p}{|y|^{n+sp}} |y|^{n+sp}|\varphi_{\rho}(y)|dy\;dx\\
+\int_{B_{\rho}} \int_{\mathbb{R}^{n}}\frac{\Big|A_{[\theta]}(t-s,y)-A_{[\theta]}(t-s,x)\Big|^p}{|x-y|^{n+sp}} |x-y|^{n+sp}\big(\varphi_\rho(x-y)+\varphi_\rho(x(s)-y)\big) dy dx\Big).
\end{eqnarray*}
But, since $x\in B_{\rho}$ and $y\in B_{\rho}$ due to the support properties of $\varphi_{\rho}$ in the first integral, and $|x-y|\le \rho, |x(s)-y|\le 2 \rho $ for the second one, we get that 
$$|y|^{n+sp}|\varphi_{\rho}(y)|\le C\rho^{n+sp}\frac{\|\varphi\|_{L^\infty}}{\rho^n}\le C\rho^{sp}\|\varphi\|_{L^\infty},$$ 
and $|x-y|^{n+sp}\big(\varphi_\rho(x-y)+\varphi_\rho(x(s)-y)\big)\le C\rho^{n+sp}\frac{\|\varphi\|_{L^\infty}}{\rho^n}\le C\rho^{sp}\|\varphi\|_{L^\infty}$. Thus:
$$\|A_{[\theta]}(t-s,\cdot)-\overline{A_{[\theta,\rho]}}(t-s,\cdot)\|_{L^p(B_\rho)}^p
\leq  C \rho^{sp} \|\varphi\|_{L^{\infty}}\int_{\mathbb{R}^{n}} \int_{\mathbb{R}^{n}}\frac{|A_{[\theta]}(t-s,x)-A_{[\theta]}(t-s,y)|^p}{|x-y|^{n+sp}}dy\;dx.
$$
With the definition of Besov spaces given in (\ref{DefinitionBesov}) we have
$$\|A_{[\theta]}(t-s,\cdot)-\overline{A_{[\theta,\rho]}}(t-s,x(s))\|_{L^p(B_\rho)}\leq C \rho^{s}\|A_{[\theta]}(t-s,\cdot)\|_{\dot{B}^{s,p}_{p}(\mathbb{R}^{n})} .$$ 
To obtain Lemma \ref{LemBesov1} it is enough to replace $s$ by $\frac{\alpha}{p}$ and to recover inequalities (\ref{EquationKeyMoins}) and (\ref{EquationKeyPlus}) we replace $s$ by $\frac{\alpha}{p}+\eta$ and by $\frac{\alpha}{p}-\eta$ respectively.
\hfill$\blacksquare$\\[5mm]
\textbf{\textit{Proof of Lemma \ref{LemBesov2}.}} The proof of this lemma follows almost the same lines than the previous one. Indeed, recalling that we now have $x\in B_{2^{k}\rho}$, we derive in the same manner that for $1<p<+\infty$ and $0<s<1$:
$$\|A_{[\theta]}(t-s,\cdot)-\overline{A_{[\theta,\rho]}}(t-s,x(s))\|_{L^p(B_{2^{k}\rho})}^p$$
\begin{eqnarray*}
&\leq & C\int_{B_{2^{k}\rho}} \int_{\mathbb{R}^{n}}\frac{|A_{[\theta]}(t-s,x)-A_{[\theta]}(t-s,y)|^p}{|x-y|^{n+sp}} |2^k\rho|^{n+sp}\frac{\|\varphi\|_{L^\infty}}{\rho^n}dy\;dx.
\end{eqnarray*}
We thus obtain
$$\|A_{[\theta]}(t-s,\cdot)-\overline{A_{[\theta,\rho]}}(t-s,\cdot)\|_{L^p(B_{2^{k}\rho})}\leq C(2^k \rho)^{s} \times 2^{\frac{kn}{p}} \times\|A_{[\theta]}(t-s,\cdot)\|_{\dot{B}^{s,p}_{p}(\mathbb{R}^{n})}.$$ 
Again, to obtain Lemma \ref{LemBesov2} it is enough to replace $s$ by $\frac{\alpha}{p}$. To recover inequalities (\ref{EquationKeyMoins1}) and (\ref{EquationKeyPlus1}) we replace $s$ by $\frac{\alpha}{p}+\eta$ and by $\frac{\alpha}{p}-\eta$ respectively. \hfill$\blacksquare$

\section*{Ackowledgments}
\hspace*{\parindent}We kindly acknowledge the two anonymous referees for their careful reading and comments which helped improving the manuscript.\\

For the second author, the article was prepared within the framework of a subsidy granted to the HSE by the Government of the Russian Federation for the implementation of the Global Competitiveness Program.


\quad\\

\begin{multicols}{2}
\begin{minipage}[r]{90mm}
Diego \textsc{Chamorro}\\[3mm]
{\footnotesize
Laboratoire de Mod\'elisation Math\'ematique d'Evry\\ 
(LaMME), UMR CNRS 8071\\ 
Universit\'e d'Evry Val d'Essonne\\[2mm]
23 Boulevard de France\\
91037 Evry Cedex\\[2mm]
diego.chamorro@univ-evry.fr
}
\end{minipage}
\begin{minipage}[r]{80mm}
St\'ephane \textsc{Menozzi}\\[3mm]
{\footnotesize
Laboratoire de Mod\'elisation Math\'ematique d'Evry\\ 
(LaMME), UMR CNRS 8071\\
Université d'Evry Val d'Essonne\\[2mm]
23 Boulevard de France\\
91037 Evry Cedex\\and \\
Laboratory of Stochastic Analysis, HSE\\ 
Moscow, Russia\\[2mm]
stephane.menozzi@univ-evry.fr
}
\end{minipage}
\end{multicols}


\begin{thebibliography}{2}

\bibitem{Adams}
D. \textsc{Adams} \& J. \textsc{Xiao}, \emph{Morrey spaces in harmonic analysis}, Ark. Mat., 50 (2012), 201–230.
\bibitem{Adams2}
D. \textsc{Adams}, \emph{Morrey Spaces}, Lecture Notes in Applied and Numerical Harmonic Analysis, Birkhäuser (2015).
\bibitem{Caffarelli}
L. \textsc{Caffarelli} \& A. \textsc{Vasseur}. \emph{Drift diffusion equations with fractional diffusion and the quasi-geostrophic equation}, Annals of Math., Vol 171 (2010), No 3, 1903-1930.


\bibitem{CCCGW}
D. \textsc{Chae} \& P. \textsc{Constantin} \& D. \textsc{C\'ordoba} \& F.
              \textsc{Gancedo} \& J. \textsc{Wu}.
     \emph{Generalized surface quasi-geostrophic equations with singular velocities},
    Comm. Pure Appl. Math., Vol 65 (2012), No 8, 1037-1066.

		
\bibitem{PGDCH}
D. \textsc{Chamorro} \& P. G. \textsc{Lemarié-Rieusset}. \emph{Quasi-geostrophic equation, nonlinear Bernstein 
inequalities and $\alpha$-stable processes}, Revista Matem\'atica Iberoamericana, Vol 28, n. 4, 1109-1122 (2012).
\bibitem{Chamorro}
D. \textsc{Chamorro}. \emph{A molecular method applied to a non-local PDE in stratified Lie groups}, J. Math. Anal. Appl. 413 (2014) 583–608.
\bibitem{DCHSM}
D. \textsc{Chamorro} \& S.  \textsc{Menozzi}. \emph{Fractional operators with singular drift: Smoothing properties and Morrey-Campanato spaces.} Revista Matem\'atica Iberoamericana. Vol 32, N°4, 1445-1499 (2016).
\bibitem{Coifmann}
R. \textsc{Coifmann} \& G. \textsc{Weiss}. \emph{Extensions of Hardy spaces and their use in analysis}, Bull Amer. Math. Soc., Vol 83, N° 4, (1977).
\bibitem{CCZV}
P. \textsc{Constantin}, M. \textsc{Coti Zelati} \& V. \textsc{Vicol}.  \emph{Uniformly attracting limit sets for the critically dissipative SQG equation}. Nonlinearity 29 (2016), N° 2, 298-318.
\bibitem{CW}
P. \textsc{Constantin} \& J. \textsc{Wu}. \emph{Regularity of Hölder continuous solutions of the supercritical quasi-geostrophic equation}, Annales de l'Institut Henri Poincaré. Analyse non linéaire. Vol 25, N°6, 1103-1110 (2008).
\bibitem{CW2}
P. \textsc{Constantin} \& J. \textsc{Wu}. \emph{H\"older continuity of solutions of supercritical dissipative hydrodynamic transport equations}, Annales de l'Institut Henri Poincaré. Analyse non linéaire. Vol 26, N°1, 159-180 (2009).
\bibitem{Cordoba}
A. \textsc{Cordoba} \& D. \textsc{Cordoba}. \emph{A maximum principle applied to quasi-geostrophic equations}, Commun. Math.
Phys. 249, 511-528  (2004).
\bibitem{Delga}
M.G. \textsc{Delgadino} \& S. \textsc{Smith}. \emph{H\"older estimates for fractional parabolic equations with critical divergence free drifts}, arXiv:1609.02691v1 (2016).
\bibitem{frie:vico:12}
S. \textsc{Friedlander} \& V. \textsc{Vicol}. \emph{ Global well-posedness for an advection-diffusion equation arising in magneto-geostrophic dynamics}, Annales de l'Institut Henri Poincar\'e (C) Analyse Non Lin\'eaire 28 (2011), no. 2, 283-301.
\bibitem{Gold}
D. \textsc{Goldberg}. \emph{A local version of real Hardy spaces}, Duke Mathematical Journal, Vol 46, N°1, (1979).
\bibitem{Jacob}
N. \textsc{Jacob}. \emph{Pseudo-Differential Operators and Markov Processes, Vol. I \& II}, Imperial College Press (2001).
\bibitem{KN}
A. \textsc{Kiselev} \& F. \textsc{Nazarov}.
\emph{A variation on a theme of Caffarelli and Vasseur}, Zapiski Nauchnykh Seminarov POMI, Vol. 370, 58–72, (2009).
\bibitem{PGLM1}
P.-G. \textsc{Lemari\'e}. \emph{Continuité sur les espaces de Besov des opérateurs définis par des intégrales singulières}. Ann. Inst. Fourier (Grenoble) 35, no. 4, 175–187 (1985).
\bibitem{Meyer}
Y. \textsc{Meyer}. \emph{Ondelettes et Opérateurs II, Opérateurs de Calder\'on Zygmund}, Hermann, (1990).
\bibitem{Peetre}
J. \textsc{Peetre}. \emph{On the theory of $L_{p,\lambda}$ spaces}. Journal of Functional Analysis, Volume 4, Issue 1, August 1969, Pages 71–87.
\bibitem{Ros}
M. \textsc{Rosenthal} \& H.-J. \textsc{Schmeisser}. \emph{On the Boundedness of Singular Integrals in Morrey Spaces and its Preduals}. JFAA, Volume 22, Issue 2, pp 462-490, (2016).
\bibitem{RS}
T.  \textsc{Runst} \& W.  \textsc{Sickel}. \emph{Sobolev spaces of fractional order, Nemytskij operators and nonlinear partial differential equations}. De Gruyter series in nonlinear analysis and applications, 3, (1996).
\bibitem{Silv}
L. \textsc{Silvestre}, V. \textsc{Vicol} \& A. \textsc{Zlatos} \textit{On the loss of continuity for super-critical drift-diffusion equations}, Archive of Rational Mechanics and Analysis 207 (2013), No 3, 845-877.
\bibitem{Stein2}
E. M. \textsc{Stein}. \emph{Harmonic Analysis}, Princeton University Press (1993).
\bibitem{Torchinski}
A. \textsc{Torchinski}. \emph{Real-Variable methods in Harmonic Analysis}, Dover, (2004).
\bibitem{Zorko}
C.T. \textsc{Zorko}.  \emph{Morrey space},  {Proc. Amer. Math. Soc.}, {98}, {586--592} (1986).

\end{thebibliography}
\end{document}